\documentclass[12pt]{book}
\usepackage{latexsym}
\usepackage{amsmath}
\usepackage{amssymb}
\usepackage{amscd}
\usepackage{amsthm}
\usepackage{xypic}
\xyoption{curve}
\usepackage{ifthen}
\usepackage{hyperref}
\usepackage{graphicx}

\theoremstyle{plain}
 \newtheorem{theorem}{Theorem}[subsection]
 \newtheorem{proposition}[theorem]{Proposition}
 \newtheorem{lemma}[theorem]{Lemma}
 \newtheorem{corollary}[theorem]{Corollary}
 
\theoremstyle{definition}
 \newtheorem{definition}[theorem]{Definition}
 \newtheorem{example}[theorem]{Example}
 
 \newtheorem{problem}[theorem]{Problem}

\theoremstyle{remark}
 \newtheorem{remark}[theorem]{Remark}

\newcommand{\RR}{\mathbb{R}}
\newcommand{\PP}{\mathbb{P}}
\newcommand{\QQ}{\mathbb{Q}}
\newcommand{\ZZ}{\mathbb{Z}}
\newcommand{\CC}{\mathbb{C}}

\newcommand{\xycenter}[1]{\begin{center} \mbox{\xymatrix{#1}}\end{center}}

\newboolean{xlabels}
\newcommand{\xlabel}[1]{
                        \label{#1}
                        \ifthenelse{\boolean{xlabels}}
                                   {\marginpar[\hfill{tiny #1}]{{\tiny #1}}}
                                    {}
                          }
 \setboolean{xlabels}{false}

\begin{document}
\frontmatter
\title{The Rationality Problem in Invariant Theory}
\author{Christian B\"ohning\\
Mathematisches Institut\\
Georg-August-Universit\"at G\"ottingen}
\date{}
\maketitle

\tableofcontents

\chapter[Introduction]{Introduction}

Invariant theory as a mathematical discipline on its own originated in England around the middle of the nineteenth century with Cayley's papers on
hyperdeterminants and his famous Memoirs on Quantics, followed by Salmon, Sylvester and Boole, and Aronhold, Clebsch and Gordan in Germany. There was also a
third school in Italy associated with the names of Brioschi, Cremona, Beltrami and Capelli. The techniques employed in this early phase, long before Hilbert
transformed the subject with his conceptual ideas, were often computational and symbolic in nature. One of the main questions was, given a linear algebraic
group $G$ and finite-dimensional $G$-representation $V$ over $\CC$, to describe the algebra of invariant polynomial functions $\CC [V]^G$ explicitly; in
fact, most attention was given classically to the case where $G=\mathrm{SL}_2 (\CC )$ or $G=\mathrm{SL}_3 (\CC )$ and $V$ is a space of binary or ternary
forms of some fixed degree.\\
Suppose now $G$ to be connected and semisimple. Today we know by work of Popov that the algebra of invariants $\CC [V]^G$ can be arbitrarily
complicated: a natural measure for its complexity is the length of its syzygy chain or in other words its homological dimension
$\mathrm{hd} (\CC[V]^G)$. Then (see e.g. \cite{Po92}, Chapter 3) it is known that if $G$ is nontrivial, then for any $n\in \mathbb{N}$, there
exists a $G$-module $V$ with
$\mathrm{hd} (\CC[V]^G) > n$ and there exist, up to isomorphism and addition of trivial direct summands, only finitely many $G$-modules with
$\mathrm{hd}(\CC[V]^G) \le n$. Moreover, the complexity of invariant rings increases quite rapidly: classically, a finite generating set and finite
set of defining relations for $\CC [\mathrm{Sym}^d (\CC^2)^{\vee }]^{\mathrm{SL}_2 (\CC )}$ was only obtained for $d\le 6$ to which the 20th
century (Dixmier \& Lazard, Shioda) contributed just $d=7,\: 8$. For $d> 8$ the homological dimension of the algebra of invariants is known to be
greater than $10$ (cf.
\cite{Po-Vi}, \S 8).\\
Thus, algebraically, one is lead to ask: when is the structure of invariants of a $G$-module $V$ as simple as possible? If we interpret this as asking when
$\CC[V]^G$ is free, i.e. has algebraically independent homogeneous generators, then, by Popov's theorem, the classification of such $V$ is a finite problem
and more information on it can be found in \cite{Po92}, \cite{Po-Vi}. One can also try to classify $G$-modules $V$ with $\CC[V]^G$ of fixed homological
dimension. However, the situation aquires a very interesting different flavour if we shift from a biregular to a birational point of view, and ask

\begin{quote}
When is $\CC (V)^G$, the field of invariant rational functions, a purely transcendental extension of $\CC$ or, as we will say, \emph{rational}?
\end{quote} 

This is the main question various aspects of which we will treat in this work. 
$\CC (V)^G$ is always of finite transcendence degree over $\CC$ (there is no Hilbert's 14th Problem birationally), and we can ask this for any linear
algebraic group $G$ whatsoever. If $G$ is not assumed to be connected, there are examples by Saltman \cite{Sa} that $\CC (V)^G$ need not even become rational
after adjunction of a number of additional indeterminates ($\CC (V)^G$ is not \emph{stably rational}). $G$ can be taken as a finite solvable group acting on
$V$ through its regular representation. This contradicts a conjecture put forward originally by Emmy Noether.\\
The quite astonishing fact, though, given the complexity of invariant rings themselves, is that no example with irrational $\CC (V)^G$ is known if $G$ is
assumed to be connected! Putting $X=V$, we can reinterpret our original question as asking: when is the quotient variety $X/G$ rational? $X/G$ is taken in
the sense of Rosenlicht and well-defined up to birational equivalence. One may replace $X$ by e.g. a rational homogeneous variety and ask the same question:
again no example of an irrational quotient $X/G$ is known if $G$ is connected.\\
The introduction of the geometric point of view is not only a reformulation, but an indispensible step for any progress on our original algebraic problem.
One may add as another example the solution to the L\"uroth Problem in dimension $2$: is an algebraic function field $L$ of transcendence degree $2$ over
$\CC$ which is contained in a purely transcendental extension of $\CC$ itself a purely transcendental extension of $\CC$? The affirmative answer follows as a
corollary of Castelnuovo's Theorem characterizing smooth projective rational surfaces as those that do not have (non-zero) holomorphic one-forms and whose
bicanonical linear system is empty. There is apparently no purely algebraic proof of this fact, though there was a time when some people tried to rewrite
the Italian birational theory of algebraic surfaces in terms of function fields. We mention that there are counter-examples to the L\"uroth Problem in
dimension $3$ and higher (cf. Artin and Mumford \cite{A-M}). So there are examples of \emph{unirational} algebraic function fields $L$ (unirational means
contained in a purely transcendental extension of $\CC$) of transcendence degree $\ge 3$ which are not themselves purely transcendental extensions of
$\CC$. There are also examples of stably rational non-rational $L$ for transcendence degree $3$ and higher \cite{Beau}, which is the solution to the
\emph{Zariski problem}. Thus we have the strict inclusions

\begin{quote}
$\{$ rational $L \}$ $\subsetneqq$ $\{$ stably rational $L \}$  $\subsetneqq$ $\{$ unirational $L \}$\, .
\end{quote}

Another reason to study quotients of the form $V/G$ (or $\PP (V) /G$) is that many moduli spaces in algebraic geometry are of this so-called \emph{linear
type}. For example, $\mathfrak{M}_g$, the moduli space of curves of genus $g$, is known to be of linear type for $1 \le g\le 6$. For example, $\mathfrak{M}_1
\simeq \PP (\mathrm{Sym}^4 (\CC^2)^{\vee}) /\mathrm{SL}_2 (\CC )$ is the ubiquitous moduli space of elliptic curves, and $\mathfrak{M}_2\simeq \PP
(\mathrm{Sym}^6 (\CC^2)^{\vee}) /\mathrm{SL}_2 (\CC )$ because a genus $2$ curve is a double cover of $\PP^1$ branched in $6$ points via its
canonical map. $\mathfrak{M}_3 \simeq
\PP (\mathrm{Sym}^4 (\CC^3)^{\vee })/\mathrm{SL}_3 (\CC )$ since a general (non-hyperelliptic) curve of genus $3$ is realized as a smooth quartic in $\PP^2$
via the canonical embedding. We do not discuss $\mathfrak{M}_4$, $\mathfrak{M}_5$, $\mathfrak{M}_6$, but just remark that certainly $\mathfrak{M}_g$ ceises
to be of linear type at some point because for $g\ge 23$, $\mathfrak{M}_g$ is not even unirational. Other examples of moduli spaces of linear type are the
moduli spaces of polarized K3 surfaces of degree $d$ for $d=2,\: 4, \: 6, \: 10$ (these classify pairs $(S, \: h)$ where $S$ is a smooth K3 surface and $h$
an ample class with $h^2 = d$ on $S$), or many moduli spaces of vector bundles. Of course, one should add to this list moduli spaces such as $\PP
(\mathrm{Sym}^d (\CC^{n+1})^{\vee })/\mathrm{SL}_{n+1} (\CC )$, the moduli space of degree $d$ hypersurfaces in $\PP^n$ (for projective equivalence) which
are of linear type by definition, and very interesting in their own right.\\
The transgression in the behaviour of $\mathfrak{M}_g$ from being rational/unirational for small $g$ and of general type for $g$ large illustrates an
important point: rational (or unirational) moduli spaces emerge as the most interesting examples (whereas the general curve of genus $g$ for large $g$ is
rather hard to put hands on as a mathematical object). In general, rational
varieties (or those close to being rational) are those that appear most frequently in applications in mathematics and make up the greatest part of one's
motivating examples in algebraic geometry, though they are only a very small portion in the class of all varieties. It is precisely the fact that they are
for the most part tangible objects and amenable to concrete study which explains their importance, and the wildness and absence of special features,
symmetries etc. which lessens the impact of the rest of varieties on the whole of mathematics.

\

What methods are there to tackle the rationality problem for $\PP(V)/G$? This is discussed in great detail in Chapters 1 and 2, so we content ourselves here
with emphasizing some general structural features and recurring problems.
\begin{itemize}
\item
If $V$ and $W$ are representations of the linear algebraic group $G$ where the generic stabilizer is trivial, then $\CC (V)^G$ and $\CC (W)^G$ are stably
equivalent, i.e. they become isomorphic after adjoining some number of indeterminates to each of them. This is the content of the so-called ''no-name lemma''
of Bogomolov and Katsylo \cite{Bo-Ka}. So the stable equivalence type is determined by the group $G$ alone in this case, and in many cases one can prove
easily that a space $V/G$ is stably rational.
\item
If one wants to prove rationality for a quotient $X/G$ ($X$ could be a linear $G$-representation or a more general $G$-variety), then, after possibly some
preparatory reduction steps consisting of taking sections for the $G$ action on $X$ and thus reducing $G$ to a smaller group and replacing $X$ by a
subvariety, virtually all the methods for proving rationality consist in introducing some fibration structure in the space $X$: one finds a
$G$-equivariant rational map $\varphi\, :\, X\dasharrow Y$ to some base variety $Y$ such that $Y/G$ is stably rational, and the generic fibre of $\varphi$ is
rational, and then one tries to use descent to prove that $X/G \dasharrow Y/G$ is birational to a Zariski bundle over $Y/G$ with rational fibre.
\item
In the examples which occur in practice where, in the situation of the previous item, $X$ is generically a $G$-vector bundle over $Y$, the map $\varphi$ can
almost always be viewed as induced from a resolution of singularities map $H\times_P F \to \mathfrak{S}$ where $\mathfrak{S}$ is a stratum of the unstable
cone in a representation $W$ of a reductive group $H \supset G$ and $F\subset W$ is some subspace which is stable under a parabolic subgroup $P \subset H$.
This method is described in detail in Section 1.3 of Chapter 1. 
\item
In the set-up of the previous two items, one almost always has to prove that the map $\varphi$ satisfies certain nondegeneracy or genericity conditions, and
this is usually a hard part of the proof. As an illustration, one can take a surface $S$ in $\PP^4$ which is the intersection of two quadric hypersurfaces
$Q_1$ and $Q_2$. To prove rationality of $S$ one projects from a line $l$ common to both $Q_1$ and $Q_2$, but one has to check that the projection is
dominant unto $\PP^2$; otherwise $S$ could be a bundle over an elliptic curve which is irrational. Checking nondegeneracy typically involves the use of
computer algebra, but special ideas are needed when e.g. one deals with an infinite number of spaces $V_n/G$, $n=1, \: 2, \dots$. A trick used in
\cite{BvB08-1} (see also Section 2.3) is to show that the data for which genericity has to be checked becomes a periodic function of $n$ over a finite field
$\mathbb{F}_p$, and then to use upper-semicontinuity over $\mathrm{Spec}(\ZZ )$ to prove nondegeneracy over $\QQ$ (or $\CC$).
\item
Finally, it would be very nice to give an example of a space $V/G$ (where $V$ as before is a linear representation of the connected linear algebraic group
$G$) which is \emph{not rational}, if such an example exists at all. A possible candidate could be given by taking $V_d$ the space of pairs of $d\times d$
matrices and $G=\mathrm{PGL}_d (\CC )$ acting on $V_d$ by simultaneous conjugation. The corresponding invariant function field $\CC (V_d)^G$ is not known to
be rational or even stably rational in general. For further information see Section 1.2.1 in Chapter 1. To determine the properties of $\CC (V_d)^G$ for
general $d$ is one of the major open and guiding problems in the subject. One should also remark that if a space $V/G$ is stably rational, then if it were
not rational, there would be practically no methods available today to prove this: the Clemens-Griffiths method of intermediate Jacobians (see \cite{Is-Pr},
Chapter 8) is limited to threefolds and the quotients $V/G$ quickly have higher dimension, the Noether-Fano-Iskovskikh-Manin method (see \emph{loc. cit.}) 
based on the study of maximal centers for birational maps has not been put to use in this context and it is hard to see how one should do it, and
Brauer-Grothendieck invariants are not sensitive to the distinction between stably rational irrational and rational varieties. 
\end{itemize}

The main known results of rationality for spaces $V/G$ can be summarized as follows: in \cite{Kat83}, \cite{Kat84}, \cite{Bogo2} and \cite{Bo-Ka} it is
proven that all quotients $\PP ( \mathrm{Sym}^d (\CC^2)^{\vee })/\mathrm{SL}_2 (\CC )$ are rational, so the problem is solved completely for binary
forms. The moduli spaces $C(d)=\PP (\mathrm{Sym}^d (\CC^3)^{\vee })/\mathrm{SL}_3 (\CC )$ of plane curves of degree $d$ are rational for $d\equiv 1$ (mod
$4$), all $d$, and for $d\equiv 1$ (mod $9$), $d\ge 19$, by \cite{Shep}, and for $d\equiv 0$ (mod $3$), $d\ge 210$, by \cite{Kat89}. This was basically
everything that was known for ternary forms prior to \cite{BvB08-1}, \cite{BvB09-1}, but there were also several rationality results for $C(d)$ for small
particular values of
$d$. Though these are somewhat sporadic, they are very valuable and should be rated rather high since rationality of $C(d)$ can be very hard to prove for
\emph{small} $d$, cf. \cite{Kat92/2}, \cite{Kat96} for the case of $C(4)$. For moduli spaces of hypersurfaces of degree $d$ in $\PP^n$ for
$n\ge 3$ much less is known: they are rational for $n=3$, $d=1, \: 2, \: 3$, and $n > 3$, $d= 1, \: 2$, which is trivial except for $n=d=3$ cf. \cite{Be}.
Likewise, spaces of mixed tensors do not seem to have been studied so far in a systematic way to my knowledge, maybe due to the smaller geometric relevance.
We should add, however, that in \cite{Shep}, the rationality of $\PP (\mathrm{Sym}^d (\CC^2)^{\vee} \otimes (\CC^2)^{\vee})/ \mathrm{SL}_2 (\CC )
\times\mathrm{SL}_2 (\CC )$, the space of pencils of binary forms of degree $d$, is proven if $d$ is even and $d\ge 10$. But for connected linear groups $G$
other than $\mathrm{SL}_n (\CC )$, $\mathrm{GL}_n (\CC )$, there are again no such good results as far as I know. The reader may turn to the surveys
\cite{Dolg1} and \cite{CT-S} for more detailed information to complement our very coarse outline.

\

We turn to the description of the contents of Chapters 1 and 2.\\
Chapter 1 introduces basic notions, gives a detailed geometric discussion of the rationality problem for the quotient of the space of pairs of $n\times n$
matrices acted on by $\mathrm{PGL}_{n} (\CC )$ through simultaneous conjugation, and presents known results for specific groups, tori, solvable groups,
special groups in the sense of S\'{e}minaire Chevalley. We then add a detailed exposition of a unifying technique for proving rationality of spaces $V/G$ that
comprises a lot of the known tricks; it uses the Hesselink stratification of the Hilbert nullcone and desingularizations of the strata in terms of
homogeneous bundles culminating in Theorem \ref{tMainTrick}. Though the method was sketched in \cite{Shep89}, it has not received such a systematic
treatment so far. In Proposition \ref{pStabRatGrassmannians} we prove a criterion for stable rationality of quotients of Grassmannians by an
$\mathrm{SL}$-action which is new, and in combination with Theorem \ref{tMainTrick} yields rationality of the moduli space of plane curves of degree $34$
(Theorem
\ref{tRationalityV34}) which was previously unknown. In section 1.4 we give a brief summary of further topics, cohomological obstructions to rationality
(unramified cohomology) and aspects of the rationality problem over fields other than $\CC$.\\ 
Chapter 2 contains a discussion of techniques available to
prove rationality of spaces
$V/G$; in part they fit in the framework of Theorem
\ref{tMainTrick}, but are presented on an elementary level with examples here which is necessary for concrete applications. In Proposition \ref{pThetaIso}
we give a reduction of the field of rational functions of the moduli space of plane curves $C$ of degree $d$ together with a theta-characteristic $\theta$
with $h^0 (C, \:
\theta )=0$ to a simpler invariant function field which is new. Finally we give an account of a method for proving rationality due to P. Katsylo, which is
based on consideration of zero loci of sections in $G$-bundles over rational homogeneous manifolds, and for which there is no good reference as far as we
know. As an application we present a proof of the rationality of the space of $7$ unordered points in $\PP^2$ modulo projectivities due to Katsylo, since the
reference is not easily accessible.\\
In Section 2.3 of Chapter 2 we summarize the results of \cite{BvB08-1}, \cite{BvB08-2}, \cite{BvB09-1} which in conclusion yield the rationality of the
moduli spaces of plane curves of degree $d$ for all but $15$ values of $d$ for which rationality remains unsettled, cf. Theorem \ref{tComprehensive}.

\

Finally I would like to thank Yuri Tschinkel for many useful discussions and for his proposal to work on the subject which turned out to be so rewarding;
furthermore I am very grateful to Fedor Bogomolov for many stimulating discussions and shaping my view of the subject. Special thanks go to Hans-Christian
Graf von Bothmer without whose mathematical and computational skills the recent results on the moduli spaces of plane curves summarized in Section 2.3 could
not have been obtained.

\mainmatter
\chapter[The rationality problem in invariant theory]{Fundamental structures in invariant theory (with an eye towards the
rationality problem)}

\section{Introduction}

In this chapter we introduce the basic notions involved in the rationality problem for invariant function fields, and discuss the action of
$\mathrm{PGL}_n (\CC )$ on pairs of $n\times n$ matrices by simultaneous conjugation as a guiding example. We give various results for specific
groups, tori, solvable groups, and special groups in the sense of S\'{e}minaire Chevalley (cf. \cite{Se58}). We introduce the Hesselink
stratification of the nullcone of a representation of a reductive group as a unifying concept for various methods for proving rationality of
quotient spaces. Together with a new criterion for the stable rationality of certain quotients of Grassmannians by an $\mathrm{SL}_n$-action
(Proposition \ref{pStabRatGrassmannians}), we obtain the rationality of the moduli space of curves of degree $34$ (Theorem \ref{tRationalityV34}).\\
In section 1.4 we give a short overview of unramified cohomology and the rationality problem over an arbitrary ground field. Apart from this section, we
work over the field of complex numbers $\CC$ throughout this text.

\section{The rationality problem}

\subsection{Quotients and fields of invariants}
Let $G$ be a linear algebraic group over $\CC$ acting (morphically) on an algebraic variety $X$.
 
\begin{definition}\xlabel{quotientdefinition}
A quotient of $X$ by the action of $G$, denoted by $X/G$, is any model of the field $\CC (X)^G$ of invariant rational functions; a quotient is
thus uniquely determined up to birational equivalence, and since we are interested in birational properties of $X/G$ here, we will also refer to
it as \emph{the} quotient of $X$ by $G$. 
\end{definition}

Note that $\CC (X)^G$ is certainly always finitely generated over $\CC$, since it is a subfield of $\CC (X)$ which is finitely generated over
$\CC$. In the context of fields there is no fourteenth problem of Hilbert (\cite{Nag}, \cite{Stein}, \cite{Muk1})!\\
Of course one would like $X/G$ to parametrize generic $G$-orbits in $X$ to be able to apply geometry.

\begin{definition}\xlabel{geometricquotient}
If $V$ is a $G$-variety, then a variety $W$ together with a morphism $\pi \, :\, V\to W$ is called a geometric quotient if 
\begin{itemize}
\item[(1)]
$\pi$ is open and surjective,
\item[(2)]
the fibres of $\pi$ are precisely the orbits of the action of $G$ on $V$,
\item[(3)]
for all open sets $U\subset W$, the map $\pi^{\ast}\, :\, \mathcal{O}_W (U) \to \mathcal{O}_V(\pi^{-1}(U))^G$ is an isomorphism.
\end{itemize}
\end{definition}

One  then has the following theorem due to Rosenlicht (\cite{Ros}, Thm. 2).

\begin{theorem} \xlabel{Rosenlicht}
There exists a nonempty $G$-stable open subset $U\subset X$ in every $G$-variety $X$ such that there is a geometric quotient for the action of $G$
on $U$.
\end{theorem}
For a modern proof, see \cite{Po-Vi} or \cite{Gross}.\\
\begin{definition}\xlabel{defrationalityversions}
\begin{itemize}
\item[(1)]
An algebraic variety $X$ is called rational if there exists a birational map $X \dasharrow \mathbb{P}^n$ for some $n$.
\item[(2)]
$X$ is called stably rational if there exists an integer $n$ such that $X\times \PP^n$ is rational.
\item[(3)]
$X$ is called unirational if there exists a dominant map $\PP^n \dasharrow X$ for some $n$.
\end{itemize}
\end{definition}
Clearly, (1)$\implies$(2)$\implies$(3) and the implications are known to be strict (\cite{A-M}, \cite{Beau}). Since (1)-(3) are properties of
the function field
$\mathbb{C}(X)$, we will also occasionally say that
$\mathbb{C}(X)$ is rational, stably rational or unirational. There are other well-known notions capturing properties of varieties which are close
to the rational varieties, notably retract rationality (\cite{Sa2}) and rational connectedness (\cite{Koll}), which we have no use for here.\\
We can now state the main problem which we are concerned with in this work.
\begin{problem}\xlabel{problemrat}
Let $G$ be a connected linear algebraic group, and let $V$ be a $G$-representation. $V$ is always assumed to be finite-dimensional.
\begin{itemize}
\item[(1)]
Is $V/G$ rational?
\item[(2)]
Is $\PP (V)/G$ rational?
\end{itemize}
\end{problem}

\begin{remark}\xlabel{remarkmainproblem}
\begin{itemize}
\item[(1)]
The existence of stably rational, non-rational varieties shows that the answer to the preceding problem is clearly no if $V$ is replaced by an
arbitrary rational variety; just take $X$ non-rational such that $X\times \CC^{\ast}$ is rational, and let the multiplicative group $\mathbb{G}_m$
act on the second factor of $X\times\CC^{\ast}$ such that $(X\times\CC^{\ast})/\mathbb{G}_m$ is birational to $X$.
\item[(2)]
By the results of Saltman (\cite{Sal}), the answer to (1) is likewise no if $G$ is not assumed to be connected; $G$ can even be taken to be a finite solvable
group acting on
$V$ via its regular representation.
\item[(3)]
The rationality of $\PP(V)/G$ implies the rationality of $V/G$. One uses the following theorem of Rosenlicht
\cite{Ros}.
\begin{theorem}\xlabel{Rosenlichtsections}
If $G$ is a connected solvable group acting on a variety $X$, then the quotient map $X \dasharrow X/G$ has a rational section $\sigma \, :\,
X/G\dasharrow X$.
\end{theorem}
In our case we have the quotient map $V/G \dasharrow \PP (V)/G$ for the action of the torus $T= \CC^{\ast}$ by homotheties on
$V/G$. If $T_{\mathrm{ineff}}$ is the ineffectivity kernel for the action of $T$ on $V/G$, the action of
$T/T_{\mathrm{ineff}}$ on $V/G$ is generically free ($T_{\mathrm{ineff}}$ coincides with the so-called stabilizer in general
position for a torus action, \cite{Po-Vi}, \S 7.2). Hence by Theorem \ref{Rosenlichtsections}, the preceding quotient map is a locally trivial
$T/T_{\mathrm{ineff}}$-principal bundle in the \emph{Zariski} topology, so that $V/G$ is birational to $\PP (V)/G \times T/T_{\mathrm{ineff}}$
which is rational if $\PP (V)/G$ is.
\item[(4)]
$\PP (V\oplus \CC ) /G$ (trivial $G$-action on $\CC$) is birational to $V/G$: Map $v\in V$ to $[(v, 1)]$ in $\PP (V\oplus\CC )$.
\end{itemize}
\end{remark}

One (and my main) motivation for Problem \ref{problemrat} comes from the fact that many moduli spaces in algebraic geometry are of the form $\PP
(V)/G$. But a solution to Problem
\ref{problemrat} or parts of it typically has diverse applications throughout algebra, representation theory and geometry. We discuss one famous
and guiding example in detail to illustrate this.

\

Let $n$ be positive integer, $G=\mathrm{GL}_n (\CC )$, and let $V=\mathfrak{gl}_n \oplus \mathfrak{gl}_n$ be two copies of the adjoint
representation of $G$ so that $V$ is the space of pairs of $n\times n$-matrices $(A, \: B)$ and $g\in G$ acts on $V$ by simultaneous conjugation:
\[
g\cdot (A,\: B) = (gAg^{-1}, \: gBg^{-1})\, .
\]
Let $K_n:=\CC (V)^G$. The question whether $K_n$ is rational is a well-known open problem. $K_n$ is known to be unirational for all $n$, stably
rational if $n$ is a divisor of $420$, rational for $n=2,\: 3, \: 4$. Excellent surveys are \cite{For02}, \cite{LeBr}. Here we just want to show
how the field $K_n$ shows up in several areas of mathematics and discuss some approaches to Problem \ref{problemrat} for $K_n$.
\begin{itemize}
\item
Let $\mathfrak{Bun}_{\mathbb{P}^2} (k,n)$ be the moduli space of stable rank $k$ vector bundles on $\PP^2$ with Chern classes $c_1=0$, $c_2=n$. It
is nonempty for $1< k\le n$. Then 
\[
\CC (\mathfrak{Bun}_{\mathbb{P}^2} (k,n) ) \simeq K_d (t_1, \dots , t_N)
\]
where the $t_i$ are new indeterminates, $d=\mathrm{gcd} (k,n)$ and $N=2nk-k^2-d^2$. In particular, for $k=n$, the field $K_n$ is the function
field of the moduli space of stable rank $n$ vector bundles on $\PP^2$ with $c_1=0$, $c_2=n$. See \cite{Kat91}.\\
The above identification arises as follows: From the monad description of vector bundles on $\PP^2$ one knows that $\CC
(\mathfrak{Bun}_{\mathbb{P}^2} (k,n) ) \simeq \CC (S_k)^G$  where $S_k$ consists of pairs $(A, B)$ of matrices such that the eigenvalues of $A$
are pairwise distinct and the rank of the commutator of $A$ and $B$ is equal to $k$. One then uses sections and the no-name lemma (see chapter 2)
to prove the above isomorphism (\cite{Kat91}).

\item
For details on the following see \cite{Pro67}, \cite{Pro76} and \cite{Pro}. Let $X=(x_{ij})$ and $Y=(y_{ij})$ be two generic $n\times n$-matrices
(the
$x_{ij}$ and
$y_{ij}$ are commuting indeterminates), and let
$R$ be the subring generated by $X$ and $Y$ inside the ring of $n\times n$ matrices with coefficients in $\CC [x_{ij}; \: y_{ij}]$. $R$ is called a
ring of generic matrices. Let $D$ be its division ring of fractions, $C$ the centre of $D$. An element in the center of $R$ is a scalar matrix
$p\cdot
\mathrm{Id}$ with $p$ a polynomial in $x_{ij}$ and $y_{ij}$ which is necessarily a polynomial invariant of pairs of matrices. $C$ is the field of
quotients of the center of $R$, thus it is a subfield of $K_n$. On the other hand, it is known that $K_n$ is generated by elements
\[
\mathrm{tr}(M_1 M_2 \dots M_{j-1}M_j)
\]
with $M_1 M_2\dots M_{j-1}M_j$ an arbitrary word in the matrices $A$ and $B$ (so each $M_i$ is either equal to $A$ or $B$). Since $D$ is a central
simple algebra of dimension $n^2$ over its centre $C$, the trace of every element of $D$ lies in $C$. Thus $K_n$ is contained in $C$, and
thus equals the centre of the generic division ring $D$.\\
\cite{Pro67} also shows that if $\bar{C}$ is the Galois extension of $C$ obtained by adjoining the roots of the characteristic polynomial of $X$
to $C$, then the Galois group is the symmetric group $\mathfrak{S}_n$, $[\bar{C} : C]=n!$ and $\bar{C}$ is a purely transcendental extension of
$\CC$. This was a stimulus to study the rationality properties of $K_n$ as a fixed field of $\mathfrak{S}_n$ acting on a rational function field
over $C$ (\cite{For79}, \cite{For80}).

\item
$K_d$ is the function field of the relative degree $g-1$ Jacobian $\mathcal{J}ac^{g-1}_d \to |\mathcal{O}_{\PP^2}(d)|_{\mathrm{smooth}}$
over the family $|\mathcal{O}_{\mathbb{P}^2}(d)|_{\mathrm{smooth}}\subset |\mathcal{O}_{\mathbb{P}^2}(d)|$ of smooth projective plane curves of
degree $d$. Here $g=(1/2)(d-1)(d-2)$ is the genus of a smooth plane curve of degree $d$. $\mathcal{J}ac^{g-1}_d$ parametrizes pairs
$(C,\:\mathcal{L})$ consisting of a smooth plane curve of degree $d$ and a line bundle $\mathcal{L}$ of degree $g-1$ on $C$. See \cite{Beau00},
section 3.

\end{itemize}
We will discuss in a little more detail now how the description of $K_d$ as the function field of a relative Jacobian over a family of plane
curves arises, and show how this can be used to give a simple geometric proof of the rationality of $K_3$ due to Michel van den Bergh.\\
The field $K_d$ is related to $\mathcal{J}ac^{g-1}_d$ via the following theorem on representations of degree $d$ plane curves as linear
determinants.

\begin{theorem}\xlabel{tlindeterminants}
Let $C$ be a smooth plane curve of degree $d$, $Jac^{g-1}_C$ its degree $g-1$ Jacobian, $\Theta \subset Jac^{g-1}_C$ the theta-divisor
corresponding to degree $g-1$ line bundles on $C$ which have nonzero global sections.\\
For each $\mathcal{L}\in Jac^{g-1}_C \backslash \Theta$, there is an exact sequence
\[
\begin{CD}
0  @>>> \mathcal{O}_{\PP^2}(-2)^d @>{A}>> \mathcal{O}_{\PP^2}(-1)^d @>>> \mathcal{L} @>>> 0
\end{CD}
\]
with $A$ a matrix of linear forms such that $\det A = F$ where $F$ is a defining equation of $C$, $I_C= (F)$.\\
Conversely, every matrix $A$ of linear forms on $\PP^2$ with $\det A =F$ gives rise to an exact sequence as before where $\mathcal{L}$ (the
cokernel of
$A$) is a line bundle on $C$ with $\mathcal{L}\in Jac^{g-1}_C \backslash \Theta$.
\end{theorem}

\begin{proof}
As is well known, the following are equivalent for a coherent sheaf on $\PP^n$ :
\begin{itemize}
\item
$\Gamma_{\ast} (\mathcal{F}):=\bigoplus_{i\in\ZZ } H^0(\PP^n, \mathcal{F}(i))$ is a Cohen-Macaulay module over the homogeneous coordinate ring $S$
of $\PP^n$.
\item
The sheaf $\mathcal{F}$ is locally Cohen-Macaulay and has trivial intermediate cohomology: $H^j (\PP^n, \mathcal{F}(t))=0$ $\forall 1\le j \le
\dim\mathrm{Supp} (\mathcal{F}) -1$, $\forall t \in\ZZ$.
\end{itemize}
Such a sheaf is called \emph{arithmetically Cohen Macaulay} (ACM).\\
Now let $\mathcal{L}$ be a degree $g-1$ line bundle on $C$ with $H^0(C, \mathcal{L})=0$. Put $\mathcal{M}:= \mathcal{L}(1)$. The ACM condition is
vacuous for line bundles on $C$ ($\mathcal{M}_x$, $x\in C$, is of course always Cohen Macaulay since $C$ is a reduced hypersurface). Moreover, one
has 
\begin{gather}\label{linearresolution}
H^0 (\mathbb{P}^2, \mathcal{M}(-1)) = H^{1} (\mathbb{P}^2, \mathcal{M}(-1)) =0 .
\end{gather}
The vanishing of $H^1$ comes from Riemann-Roch which yields $\chi (\mathcal{L})= \chi (\mathcal{M}(-1))=0$.\\
Since $\mathcal{M}$ is ACM, $\dim \mathrm{Supp} (\mathcal{M} ) + \mathrm{proj.dim} \mathcal{M} = \dim\mathbb{P}^2$ by the Auslander-Buchsbaum
formula, whence by Hilbert's syzygy theorem, $\mathcal{M}$ has a minimal graded free resolution
\[
\begin{CD}
0 \to \bigoplus_{i=1}^r \mathcal{O} (-f_i) @>{A}>> \bigoplus_{i=1}^r \mathcal{O}(-e_i) @>>> \mathcal{M} \to 0
\end{CD}
\]
where, moreover, one has $e_i\ge 0$ for all $i$ since $H^0 (\mathbb{P}^2, \mathcal{M}(-1))=0$. The support of $\mathcal{M}$, the curve $C$, is
defined by $\det A =0$ set-theoretically, whence $\det A$ is a power of $F$. If one localizes $A$ at the generic point of $C$, the above
exact sequence together with the structure theorem of matrices over a principal ideal domain yields $\det A = F$.\\
Now condition \ref{linearresolution} yields that one has $e_i=0$, all $i$, and $f_j=1$, all $j$. Namely, the condition $H^1(C, \mathcal{M}(-1))=0$
means that $\mathcal{M}$ is a $0$-regular sheaf in the sense of Castelnuovo-Mumford regularity (\cite{Mum2}) whence
\begin{quote}
$\mathcal{M}$ is spanned by $H^0(\mathcal{M})$ and for all $j\ge 0$
\[
H^0(\mathbb{P}^2, \mathcal{O}(1)) \otimes H^0(\PP^2, \mathcal{M}(j)) \to H^0(\PP^2, \mathcal{M}(j+1))
\]
is surjective.
\end{quote}
Thus $\mathcal{M}$ has the minimal graded free resolution
\[
\begin{CD}
0 \to \bigoplus_{i=1}^r \mathcal{O} (-f_i) @>{A}>> \bigoplus_{i=1}^r \mathcal{O} @>>> \mathcal{M} \to 0
\end{CD}
\]
with $r=h^0 (\mathcal{M})$. We get the exact sequence
\[
\begin{CD}
H^1 (\PP^2, \mathcal{M}(-1)) @>>> \bigoplus_{i=1}^r H^0 (\PP^2, \mathcal{O}(f_i-2)) @>>> H^0 (\PP^2, \mathcal{O}(-2))^r
\end{CD}
\]
which together with the fact that the $f_i$ must be positive (the map induced by $\bigoplus_{i=1}^r \mathcal{O} \to \mathcal{M}$ on $H^0$ is an
isomorphism) implies that we must have $f_i=1$ for all $i$. It also follows that $d=r$ since $\det A = F$.

\

Conversely, suppose that $A$ is a $d$ by $d$ matrix of linear forms on $\PP^2$ with $\det A = F$ where $F$ is a defining equation of the smooth
curve $C$. Then one has an exact sequence
\[
\begin{CD}
0 @>>> \mathcal{O}(-2)^d @>>> \mathcal{O}(-1)^d @>>> \mathcal{L} @>>> 0
\end{CD}
\]
where $\mathcal{L}$ is an ACM sheaf on $C$ of rank $1$, thus a line bundle. By the exact sequence, $H^0(\mathcal{L})= H^1(\mathcal{L})=0$ whence
$\deg \mathcal{L}= g-1$, by Riemann-Roch.
\end{proof}

\begin{corollary}\xlabel{clindeterminants}
Let $V_d$ be the vector space of $d$ by $d$ matrices of linear forms on $\PP^2$. $\mathrm{GL}_d (\CC) \times \mathrm{GL}_d (\CC)$ acts on $V_d$ by
$(M_1, M_2)\cdot A:= M_1 A M_2^{-1}$. Then $V_d / \mathrm{GL}_d (\CC) \times \mathrm{GL}_d (\CC)$ is birational to $\mathcal{J}ac^{g-1}_d$, the
relative Jacobian of degree $g-1$ line bundles over the space of smooth degree $d$ curves $C$ in $\PP^2$.
\end{corollary}

Writing $A$ in $V_d$ as 
\[
A=A_0 x_0 + A_1 x_1 + A_2 x_2, \quad A_i\in \mathrm{Mat}_{d\times d}(\CC),
\] 
we may identify $V_d$ with the space of triples $(A_0,\: A_1, \: A_2)$ of $d\times d$ scalar matrices $A_i$ where $(M_1, M_2) \in \mathrm{GL}_d
(\CC) \times \mathrm{GL}_d (\CC)$ acts as
\[
(M_1, M_2) \cdot (A_0, \: A_1,\: A_2) = (M_1 A_0 M_2^{-1}, \: M_1 A_1 M_2^{-1}, \: M_1 A_2 M_2^{-1}) \, .
\]
The subvariety
\[
\left\{ (\mathrm{Id}, \: B_1, \: B_2) \, |\, B_1, \: B_2 \in  \mathrm{Mat}_{d\times d}(\CC) \right\} \subset \mathrm{Mat}_{d\times d}(\CC) \times
\mathrm{Mat}_{d\times d}(\CC) \times \mathrm{Mat}_{d\times d}(\CC)
\]
is a $(\mathrm{GL}_{d} (\CC ) \times \mathrm{GL}_d (\CC ) , \: \mathrm{GL}_d (\CC ))$-section in the sense of Chapter 2, 2.2.1. Hence $K_d$, the
field of invariants for the action of $\mathrm{GL}_d (\CC)$ by simultaneous conjugation on pairs of matrices, is the function field of
$\mathcal{J}ac^{g-1}_d$.

\begin{remark} \xlabel{rotherJacobians}
Instead of $\mathcal{J}ac^{g-1}_d$ it is occasionally useful to work with other relative Jacobians with the same function field $K_d$: since a
line in $\PP^2$ cuts out a divisor of degree $d$ on a smooth plane curve $C$ of degree $d$, we have for $d$ odd 
\[
\mathcal{J}ac^{g-1}_d \simeq \mathcal{J}ac^{0}_d
\]
since $g-1 = \frac{1}{2}d (d-3)$; also, in general,
\[
\mathcal{J}ac^{g-1}_d \simeq \mathcal{J}ac^{ d \choose 2 }_d\, .
\]
Also note that $\mathcal{J}ac^{g}_d$ is rational since it is birational to a (birationally trivial) projective bundle over $\mathrm{Sym}^{g}\PP^2$
which  is rational, but this yields no conclusion for $\mathcal{J}ac^{g-1}_d$. 
\end{remark}

\begin{theorem}\xlabel{tratK3}
The field $K_3$, i.e. the function field of $\mathcal{J}ac^0_3$, is rational.
\end{theorem}
\begin{proof}
We follow \cite{vdBer}. We have to prove that the variety $\mathcal{J}ac^0_3$, parametrizing pairs $(C, \mathcal{L})$, where $C$ is a smooth plane
cubic and $\mathcal{L}$ is a line bundle of degree $0$ on $C$, is rational. Fix once and for all a line $l\subset\PP^2$. Let $\mathcal{L}$ be
represented by a divisor $D$ of degree $0$ on $C$. For a general curve $C$, $l$ intersects $C$ in three points $p_1,\: p_2,\: p_3$ (uniquely
defined by $C$ \emph{up to order}), and since by Riemann Roch $h^0 (C, \mathcal{O}(p_i+D))=1$, there are uniquely determined points $q_1, \: q_2,
\: q_3$ on
$C$ with
$p_i +D
\equiv q_i$, the symbol $\equiv$ denoting linear equivalence. Thus 
\[
q_i + p_j +D \equiv q_j + p_i + D \implies q_i + p_j \equiv q_j + p_i \; \forall i,\: j \in \{ 1,\: 2, \: 3\} \, ,
\]
and if $r_{ij}$ denotes the third point of intersection of the line $\overline{p_iq_j}$ through $p_i$ and $q_j$ with $C$, clearly
\[
p_i + q_j + r_{ij} \equiv p_j + q_i + r_{ji}\, ,
\]
whence $r_{ij} \equiv r_{ji}$, and by Riemann Roch in fact $r_{ij} = r_{ji}$. Thus $r_{ij} = \overline{p_iq_j} \cap \overline{p_jq_i}$ lies on
$C$, and we get nine points: $p_1,\: p_2, \: p_3$, $q_1, \: q_2, \: q_3$, and $r_{12}, \: r_{13}, \: r_{23}$.\\
Conversely, given three arbitrary points $p_1, \: p_2, \: p_3$ on $l$, and three further points $q_1$, $q_2$, $q_3$ in $\mathbb{P}^2$, we may set
for $i < j$ $r_{ij}:= \overline{p_iq_j} \cap \overline{p_jq_i}$ and find a cubic curve $C$ through all of the $p_i$, $q_j$, $r_{ij}$. Then $D:=
q_1 - p_1$ is a degree $0$ divisor on $C$. Applying the preceding construction, we get back the points we started with.\\
However note that the cubic curve $C$ through $p_i$, $q_j$, $r_{ij}$ as before is not unique: generally, there is a whole $\PP^1$ of such curves
$C$. This is because any cubic $C$ passing through the eight points $p_1, \: p_2, \: p_3, \: q_1, \: q_2, \: q_3, \: r_{12}, \: r_{13}$ passes
through the ninth point $r_{23}$ since 
\begin{gather*}
q_1 + p_2 + r_{12} \equiv p_1 + q_2 +r_{12} \quad \mathrm{and} \quad q_1 + p_3 + r_{13} \equiv p_1 + q_3 +r_{13}\\
\mathrm{implies} \quad p_2 + q_3 \equiv p_3 + q_2 \quad \mathrm{on}\; C 
\end{gather*}
whence the intersection point $r_{23} = \overline{p_2q_3} \cap \overline{p_3q_2}$ necessarily lies on $C$. By explicit computation one may check
that for generic choice of the $p_i$ and $q_j$ one gets indeed a pencil of cubic curves $C$.\\
The above can be summarized as follows:
\begin{itemize}
\item
Let $\mathcal{P}$ be the parameter space inside $l^3 \times (\PP^2)^3 \times \PP H^0 (\PP^2, \mathcal{O}(3))$ consisting
of triples $((p_1,\: p_2,\: p_3), \; (q_1, \: q_2,\: q_3), \; C)$ where the $p_i$ are three points on $l$, the $q_j$ are three further points in
$\PP^2$, and $C$ is a cubic curve through the $p_i$, $q_j$, and $r_{ij}:=\overline{p_iq_j} \cap \overline{p_jq_i}$. Let the symmetric group
$\mathfrak{S}_3$ act on $\mathcal{P}$ via
\begin{gather*}
\sigma\cdot ((p_1,\: p_2,\: p_3), \; (q_1, \: q_2,\: q_3), \; C)\\
 := ((p_{\sigma (1)},\: p_{\sigma (2)},\: p_{\sigma (3)}), \; (q_{\sigma (1)}, \:
q_{\sigma (2)},\: q_{\sigma (3)}), \; C)\, .
\end{gather*}
Then $\mathcal{J}ac^0_3$ is birational to $\mathcal{P}/\mathfrak{S}_3$.
\item
Let $\mathcal{Q}$ be the parameter space $l^3 \times (\PP^2)^3$ of three points $p_i$ on $l$ and three additional points
$q_i$ in $\PP^2$. $\mathfrak{S}_3$ acts on $\mathcal{Q}$:
\[
\sigma\cdot ((p_1,\: p_2,\: p_3), \; (q_1, \: q_2,\: q_3)) := ((p_{\sigma (1)},\: p_{\sigma (2)},\: p_{\sigma (3)}), \; (q_{\sigma (1)}, \:
q_{\sigma (2)},\: q_{\sigma (3)})\, .
\]
We have a forgetful map $\mathcal{P} \to \mathcal{Q}$ (the field of rational functions $\CC (\mathcal{P})$ is a purely transcendental extension of $\CC (\mathcal{Q})$
given by adjoining the solutions of a set of linear equations). Passing to the quotients, we get a map $\pi : \mathcal{P} /\mathfrak{S}_3 \dasharrow
\mathcal{Q}/\mathfrak{S}_3$. Since the action of
$\mathfrak{S}_3$ on
$\mathcal{Q}$ is generically free,
$\pi$ is generically a
$\mathbb{P}^1$-bundle (in the classical topology, i.e. a conic bundle).
\end{itemize}
To conclude the proof, it suffices to remark that the conic bundle $\pi$ has a rational section, hence is birationally trivial: indeed, one just
has to fix one further point $x \in \PP^2$, and assigns to points $((p_1,\: p_2,\: p_3), \; (q_1, \: q_2,\: q_3))$ the triple $((p_1,\: p_2,\:
p_3), \; (q_1, \: q_2,\: q_3), \; C)$ where $C$ is the unique cubic passing through $p_i$, $q_j$, $r_{ij}:=\overline{p_iq_j} \cap
\overline{p_jq_i}$ and the point $x$. Moreover, the base $\mathcal{Q}/\mathfrak{S}_3$ is clearly rational: since the action of
$\mathfrak{S}_3$ on $\CC^3\oplus \CC^3\oplus \CC^3$ by permuting the factors is generically free, one sees from the no-name lemma (cf.
Chapter 2, subsection 2.2.2) and the existence of sections for torus actions \ref{Rosenlichtsections}, that $\CC(l^3 \times
(\PP^2)^3)^{\mathfrak{S}_3}$ is a purely transcendental extension of $\CC ((\PP^2)^3)^{\mathfrak{S}_3}$ which is rational.
\end{proof}

This proof is more geometric (and easier from my point of view) than the one given in \cite{For79}. Since the projective geometry of plane
quartics is quite rich, we would like to ask whether one can also obtain the result of \cite{For80} in this way.

\begin{problem} \xlabel{pFormanekquartics}
Can one prove the rationality of $K_4$ using its identification with $\mathcal{J}ac^2_4$ and the classical projective geometry of plane quartics?
\end{problem}

The following remark shows that the stable rationality of $K_d$ does not follow from a straightforward argument that is close at hand.

\begin{remark}\xlabel{rstabratnoteasy}
Put $k= {d\choose 2}$ and in $\mathrm{Sym}^{k}\PP^2 \times \mathbb{P}(H^0 (\PP^2, \mathcal{O}(d)))$ consider the incidence correspondence $X$ given
by the rule that $k$ unordered points in $\mathrm{Sym}^k\PP^2$ lie on a plane curve of degree $d$ in $\PP (H^0 (\PP^2, \mathcal{O}(d)))$. $X$ is
generically the projectivisation of a vector bundle over $\mathrm{Sym}^k\PP^2$ hence rational. On the other hand, one has also the natural map
$X\to \mathcal{J}ac^{k}_d$, assigning to a pair $(D, \: C)$ the point $( | D| , \: C)$ in $\mathcal{J}ac^{k}_d$, which makes $X$ a $\PP^N$-bundle
in the \emph{classical} or \emph{\'{e}tale} topology over the dense open subset $U \subset \mathcal{J}ac^{k}_d$ consisting of pairs $(| D|, \: C)$
with $|D|$ a non-special divisor class.\\
However, unfortunately, this is not a projective bundle in the \emph{Zariski} topology: if, to the contrary, this was the case, let $\sigma \, :\,
\mathcal{J}ac^k_d \dasharrow X$ be a rational section. Then, if $u\, :\, U_d \to |\mathcal{O}(d)|_{\mathrm{smooth}}$ is the universal curve, the
pull-back of $(|D| , \sigma (|D| ))$ on $\mathcal{J}ac^k_d \times_{|\mathcal{O}(d)|_{\mathrm{smooth}}} X$ to $\mathcal{J}ac^k_d
\times_{|\mathcal{O}(d)|_{\mathrm{smooth}}} U_d$ would give a universal divisor or Poincar\'{e} line bundle on $\Omega \times_V u^{-1}(V)$ where
$V\subset |\mathcal{O}(d)|_{\mathrm{smooth}}$ is some dense open set, and $\Omega$ some dense open set in $\mathcal{J}ac^k_d|_V$.
But by results of Mestrano and Ramanan (\cite{Me-Ra}, Lemma 2.1 and Corollary 2.8), 
\begin{quote}
a Poincare bundle on $\mathcal{J}ac^{r}_d \times_{|\mathcal{O}(d)|_{\mathrm{smooth}}} U_d$ (or on $\Omega
\times_{|\mathcal{O}(d)|_{\mathrm{smooth}}} u^{-1}(V)$ as above) exists if and only if $1-g+r$
and
$d$ are coprime, where $g=\frac{1}{2}(d-1)(d-2)$ is the genus.
\end{quote}
Hence in all the cases we are interested in, $X \to \mathcal{J}ac^{k}_d$ is generically a nontrivial Severi-Brauer scheme.
\end{remark}

\subsection{Results for specific groups}
Here we collect some results in the direction of Problem \ref{problemrat} which express exclusively properties of the group $G$ acting, for
specific groups $G$, and are independent of the particular $G$-representation $V$.\\
The following theorem is due to Miyata \cite{Mi}. 
\begin{theorem}\xlabel{Miyata}
Let $G$ be isomorphic to a subgroup of the Borel group $B_n\subset \mathrm{GL}_n\, (\CC )$ of invertible upper triangular
matrices. Then the field of invariant rational functions $\CC (V)^G$ for the $G$-module $V=\CC^n$ is a purely
transcendental extension of $\CC$.
\end{theorem}
\begin{proof}
The proof is an immediate application of the following
\begin{quote}
\bf{Claim}.\mdseries If $k$ is a field and $G$ a group of automorphisms of the polynomial ring $k[t]$ in one indeterminate
which transforms
$k$ into itself , then there is an invariant $p\in k[t]^G$ such that $k(t)^G=k^G (p)$.
\end{quote}
To prove the claim, note that $k(t)^G$ is the field of fractions of $k[t]^G$: write $f\in k(t)^G$ as $f=u/v$, $u,\: v\in
k[t]$ without common factor. After passing to the reciprocal if necessary, we may assume $\mathrm{deg}(u)\ge \deg (v) >0$
and apply the division algorithm in $k[t]$ to write
\[
u=qv + r\, ,
\]
$q,\: r\in k[t]$, $\deg (r) < \deg (v)$ whence $q$ and $r$ are uniquely determined by these requirements. Since $f$ is
invariant, $G$ acts on both $u$ and $v$ via a certain character $\chi\, :\, G\to \CC^{\ast}$, and the uniqueness of $r$
and $q$ implies that $r$ is a weight vector of $G$ for the character $\chi$, and $q$ is an (absolute) $G$-invariant. Since
\[
\frac{u}{v} = q + \frac{r}{v} , \; \frac{r}{v}\in k(t)^G ,
\]
and $\deg (r) + \deg (v) < \deg (u) +\deg (v)$, one obtains the statement by induction on $\deg (u) + \deg (v)$, the case
$\deg (u) +\deg (v)=0$ being trivial.\\
Now if $k[t]^G\subset k$, the claim is obvious. Otherwise, we take $p\in k[t]^G\backslash k$ of minimal degree. Then if
$f$ is in $k[t]^G$, one writes $f=pq+r$ with $\deg (r) < \deg (p)$ as before, and by uniqueness of quotient and
remainder, $q$ and $r$ are $G$-invariant polynomials in $k[t]$. Thus, by the choice of $p$, $r\in k^G$ and $\deg (q) <
\deg (f)$. Again by induction on the degree of $f$ we obtain $f\in k^G[p]$. This means $k[t]^G=k^G[p]$, and since we have
seen that $k(t)^G$ is the field of fractions of $k[t]^G$ the assertion of the claim follows.\\
To prove the theorem, one applies the claim to $k= \mathbb{C}(x_1,\dots , x_{n-1})$, $t=x_n$ where $x_1, \dots , x_n$ are
coordinates on $V=\CC^n$, and concludes by induction on the number of variables.
\end{proof}

\begin{corollary}\xlabel{Miyatacorollary}
If $V$ is a finite dimensional linear representation of either
\begin{itemize}
\item
an abelian group $G\subset \mathrm{GL}(V)$ consisting of semi-simple elements (e.g. if $G$ is finite)
\item
or a connected solvable group $G$,
\end{itemize}
then $\mathbb{C}(V)^G$ is a purely transcendental extension of $\CC$.
\end{corollary}
\begin{proof}
Simultaneous diagonalizability of commuting semisimple elements, or Lie-Kolchin theorem, respectively.
\end{proof}

\begin{remark}\xlabel{Miyatanonclosed}
Note that the statement and proof of Theorem \ref{Miyata} remain valid if one works, instead of over $\CC$, over a possibly nonclosed ground
field; the Corollary \ref{Miyatacorollary} becomes false in general, however, because one needs the algebraic closedness to make the actions
triangular: for example, let $G$ be a cyclic group of order $p=47$, and let $G$ act on $\QQ (x_1, \dots , x_p)$ by permuting the variables cyclically. It is
known (cf. \cite{Swan}) that $\QQ (x_1, \dots , x_p)^G$ is not rational \emph{over} $\QQ$.
\end{remark}

For semi-simple groups the only truely complete rationality result is the following due to P. Katsylo and F. Bogomolov.

\begin{theorem}\xlabel{tKatsyloSL2}
The moduli spaces $\PP (\mathrm{Sym}^d \CC^2)/ \mathrm{SL}_2 (\CC )$ of $d$ unordered points in $\PP^1$ are rational for all $d$.
\end{theorem}

See \cite{Kat84}, \cite{Bo-Ka}, \cite{Bogo2} for a proof. It should be mentioned that this result is also used in the recent work of Kim and Pandharipande
\cite{Ki-Pa}. There they prove the rationality of the moduli space
\[
\overline{M}_{0,\: n} (X, \: \beta ) 
\]
of $n$-pointed genus $0$ stable maps of class $\beta\in H_2 (X, \: \ZZ )$ into a rational homogeneous variety $X=G/P$. Here $\overline{M}_{g, \: n} (X,\: 
\beta )$ parametrizes data
\[
[\mu \, :\, C \to X ; \: p_1, \dots , p_n ]
\] 
where $C$ is a complex, projective, connected, reduced, nodal curve of arithmetic genus $g$, $p_1, \dots , p_n$ are distinct points in the smooth locus of
$C$, the map $\mu$ has no infinitesimal automorphisms and $\mu_{\ast }[C] = \beta$.\\
Despite the very small number of general rationality results, one has some satisfactory
information with regard to the important question of existence of rational sections. We have already seen Rosenlicht's theorem \ref{Rosenlichtsections}. Let
us recall quickly the theory of \emph{special groups} cf. \cite{Se56}, \cite{Se58}, \cite{Groth58}.

\begin{definition}\xlabel{dprincbundle}
Let $G$ be an algebraic group, and $\pi \, :\, P \to X$ a morphism of algebraic varieties. Let $P$ be equipped with a \emph{right} $G$-action and
suppose
$\pi$ is constant on
$G$ orbits. Then $P$ is called \emph{a $G$-principal bundle in the \'{e}tale topology} (or \emph{locally isotrivial fibre space with typical fibre}
 $G$ or $G$\emph{-torsor}) if for every point $x\in X$ there is a Zariski open neighborhood $U\ni x$ and an \'{e}tale cover $f\, :\, U' \to U$
such that the pull-back $f^{\ast}(P|_U) \to U'$ is $G$-isomorphic to the trivial fibering $U'\times G \to U'$.\\
$P$ is called \emph{a $G$-principal bundle in the Zariski topology} if furthermore every $x\in X$ has a Zariski open neighborhood $U$ such that
$P|_U$ is trivial.
\end{definition}

\begin{definition}\xlabel{dspecialgroups}
An algebraic group $G$ is called \emph{special} if every $G$-principal bundle in the \'{e}tale topology is Zariski locally trivial.
\end{definition}

The main results of interest to us are summarized in the following

\begin{theorem}\xlabel{tspecialgroups}
\begin{itemize}
\item[(a)]
The general linear group $\mathrm{GL}_n (\CC )$ is special.
\item[(b)]
A closed subgroup $G\subset \mathrm{GL}_n (\CC )$ is special if and only if the quotient map $\mathrm{GL}_n (\CC ) \to \mathrm{GL}_n (\CC ) /G$ is
a Zariski locally trivial $G$-principal bundle (equivalently, if and only if $\mathrm{GL}_n (\CC ) \to \mathrm{GL}_n (\CC ) /G$ has a rational
section).
\item[(c)]
The groups $\mathrm{SL}_n (\CC )$, $\mathrm{Sp}_n (\CC )$ and all connected linear solvable groups are special.
\item[(d)]
If $G$ is a linear algebraic group, $H$ a normal subgroup, and if $H$ and $G/H$ are special, then $G$ is special. In particular, any connected
linear algebraic group the semisimple part of which is a direct product of groups of types $\mathrm{SL}$ or $\mathrm{Sp}$ is special.
\end{itemize}
\end{theorem}

\begin{proof}
(a):The proof uses the method of taking averages of group cocycles.\\
Let $P \to X$ be a $G$-principal bundle, and let $X' \to X$ be a finite \'{e}tale cover, which we can assume to be Galois with Galois group $\Gamma$,
such that
$P$ becomes trivial on $X'$:
\[
\begin{CD}
X'\times G @>>> P\\
@V{\pi'}VV      @V{\pi}VV\\
X' @>f>>    X
\end{CD}
\]
We recall the bijective correspondence between the set of isomorphism classes of $G$-principal bundles (in the \'{e}tale topology) on $X$ which
become trivial when pulled-back to $X'$ and the elements of the nonabelian group cohomology set $H^1 (\Gamma , \mathrm{Mor} (X', \, G))$ with
marked point. Here $\mathrm{Mor}(X', \, G)$ is the group of morphisms of $X'$ into $G$; and $\Gamma$, which we assume to operate on the right on
$X'$, acts on $\mathrm{Mor}(X', \, G)$ (on the left) by
\[
(\sigma \cdot \varphi ) (x') := \varphi (x' \cdot \sigma )\, .
\]
Elements of $H^1 (\Gamma , \mathrm{Mor} (X', \, G))$ are by definition $1$-cocycles of $\Gamma$ with values in $\mathrm{Mor} (X', \, G)$ modulo an
equivalence relation; a $1$-cocycle is a map $\sigma \mapsto \varphi_{\sigma}$ from $\Gamma$ to $\mathrm{Mor}(X', \, G)$ satisfying
\[
\varphi_{\sigma \tau} = (\varphi_{\tau })^{\sigma } \varphi_{\sigma }
\]
where $(\cdot )^{\sigma}$ denotes the action of $\sigma$. Two 1-cocycles $( \varphi_{\sigma})$, $(\varphi_{\sigma}')$ are cohomologous if there is
an $a \in \mathrm{Mor}(X', \, G)$ such that 
\[
\varphi_{\sigma }' = a^{\sigma } \varphi_{\sigma} a^{-1}, \quad \mathrm{all} \; \sigma \, .
\]
Now in the above pull-back diagram $X'\times G$ is a Galois cover of $P$ with Galois group $\Gamma$, and $P = (X'\times G)/\Gamma$. $\Gamma$ acts
on $X'\times G$ compatibly with the projection $\pi'$ to $X'$ and the operation of $G$ whence
\[
(x', g) \cdot \sigma = (x' \cdot \sigma, \: \varphi_{\sigma }(x') \cdot g)
\]
and the associativity gives the required cocycle condition for $(\varphi_{\sigma })$. Conversely, the datum of a $1$-cocycle $(\varphi_{\sigma })$
determines an operation of $\Gamma$ on $X'\times G$ and one may define $P$ on $X$ as the quotient. The condition that two $1$-cocycles are
cohomologous means precisely that the $G$-principal bundles so obtained are isomorphic. Note that the condition that $(\varphi_{\sigma})$ and
$(\varphi_{\sigma}')$ are cohomologous means precisely that the isomorphism $(x', \: g) \mapsto (x', \: a(x') g)$ between trivial $G$-principal
bundles on $X'$ descends to the $G$-principal bundles on $X$ defined by $(\varphi_{\sigma})$ and $(\varphi_{\sigma}')$ on $X$, respectively.\\ We
turn to the proof of (a) of Theorem
\ref{tspecialgroups}. Thus let
$P
\to X$ be a
$\mathrm{GL}_n (\CC )$-principal bundle in the
\'{e}tale topology, and let $U\ni x$ be an open neighborhood, $f\, :\, U' \to U$ a Galois cover with group $\Gamma$ on which $P$ is trivial. The
above considerations show that, if $\mathcal{O}_{f^{-1}(x)}$ is the semi-local ring of the fibre over $x$, then the set of isomorphism classes of
$\mathrm{GL}_n (\CC )$-principal bundles on a Zariski neighborhood of $x$ which become trivial on a Zariski neighborhood of $f^{-1}(x)$ are
identified with the cohomology set $H^1 (\Gamma , \: \mathrm{GL}_n (\mathcal{O}_{f^{-1}(x)}))$. Let $x'$ be a point of the fibre $f^{-1}(x)$ and
choose a matrix $b\in \mathrm{Mat}_{n\times n} (\mathcal{O}_{f^{-1}(x)})$ which is the identity in $x'$ and the zero matrix in the other points of
$f^{-1}(x)$. If $(\varphi_{\sigma })$ is a $1$-cocycle representing the germ of $P$ in $x$ one puts
\[
a = \sum_{\tau \in \Gamma } \tau (b) \varphi_{\tau} \, .
\]
By definition, this is invertible in each point of the fibre $f^{-1}(x)$, thus belongs to $\mathrm{GL}_n (\mathcal{O}_{f^{-1}(x)} )$. Since
\[
a^{\sigma} \varphi_{\sigma } = \sum_{\tau \in \Gamma} \sigma (\tau (b))  (\varphi_{\tau })^{\sigma } \varphi_{\sigma } = \sum_{\tau \in\Gamma}
( \sigma\tau ) (b) \varphi_{\sigma\tau} =a \, ,
\]
we have $a^{\sigma} \varphi_{\sigma } a^{-1} = 1$, so our $\mathrm{GL}_n (\CC )$-principal bundle is trivial in a Zariski neighborhood of $x$.\\
Using the correspondence between $\mathrm{GL}_n (\CC )$-principal bundles in the \'{e}tale resp. Zariski topology and vector bundles in the
\'{e}tale resp. Zariski topology (given by passing to the associated fibre bundles with typical fibre $\CC^n$, and conversely associated frame
bundles), we obtain the fact which is fundamental to many techniques for proving rationality, that every vector bundle in the \'{e}tale topology
is a vector bundle in the Zariski topology.

\

(b): The proof consists in the trick of extension and reduction of the structure group.\\
By definition, if $G\subset \mathrm{GL}_n (\CC )$ is
special, then
$\mathrm{GL}_n (\CC )
\to
\mathrm{GL}_n (\CC ) /G$ is Zariski locally trivial. Conversely, suppose that $\mathrm{GL}_n (\CC ) \to \mathrm{GL}_n (\CC ) /G$ is Zariski
locally trivial, and let
$P\to X$ be a $G$-principal bundle in the \'{e}tale topology. We have an associated fibre bundle $Q:=P \times^G \mathrm{GL}_n (\CC )$ which is a
$\mathrm{GL}_n (\CC )$-principal bundle \emph{in the Zariski topology} by part (a).\\
Now $Q$ is a $G$-principal bundle in the Zariski topology over $Q\times^{\mathrm{GL}_n (\CC )} (\mathrm{GL}_n (\CC )/G)$ because $Q$ itself is
Zariski locally trivial and $\mathrm{GL}_n (\CC ) \to \mathrm{GL}_n (\CC ) /G$ has a rational section by assumption. Now 
\[
Q\times^{\mathrm{GL}_n (\CC )} (\mathrm{GL}_n (\CC )/G) = P \times^G (\mathrm{GL}_n (\CC )/G)
\] 
has a canonical section $\sigma \, :\, X \to P \times^G (\mathrm{GL}_n (\CC )/G)$ since $G$ leaves the coset corresponding to the identity in
$\mathrm{GL}_n (\CC )/G$ invariant. Then $P$ is the pull-back of $Q \to P \times^G (\mathrm{GL}_n (\CC )/G)$ via $\sigma$:
\[
\begin{CD}
P  @>>>  Q = P \times^G \mathrm{GL}_n (\CC )\\
@VVV      @VVV\\
X  @>{\sigma }>>  P \times^G (\mathrm{GL}_n (\CC )/G)
\end{CD}
\]
Thus the fact that $P \to X$ is Zariski locally trivial follows from the fact that $Q  \to P \times^G (\mathrm{GL}_n (\CC )/G)$ has this property. 

\

(c): For connected linear solvable groups, this follows from part (d) to be proven below since a connected solvable group is a successive
extension of groups of type $\mathbb{G}_m$ and $\mathbb{G}_a$. Remark that both $\mathbb{G}_m$ and $\mathbb{G}_a$ are special since
$\mathbb{G}_m=\mathrm{GL}_1 (\CC )$ and the natural map $\mathrm{GL}_2 (\CC ) \to \mathrm{GL}_2 (\CC )/\mathbb{G}_a$ has a rational section. The projection
$\mathrm{GL}_n (\CC )
\to
\mathrm{GL}_n (\CC )/\mathrm{SL}_n (\CC )$ has a section given by assigning to a coset $g\mathrm{SL}_n (\CC )$ the matrix
\[
\mathrm{diag} (\det g, \: 1, \dots , 1)\, .
\]
Finally, the projection $\mathrm{GL}_n (\CC ) \to \mathrm{GL}_n (\CC )/\mathrm{Sp}_n (\CC )$, $n=2m$, has a section since $\mathrm{GL}_n (\CC
)/\mathrm{Sp}_n (\CC )$ is the space of nondegenerate skew-symmetric bilinear forms on $\CC^{2m}$ and the generic
skew-symmetric form
\[
\sum_{1\le i < j \le n} t_{ij} (x_iy_j - y_i x_j) 
\]
with indeterminate coefficients $t_{ij}$ can be reduced to the canonical form $\sum_{k=1}^m (x_{2k-1}y_{2k}-y_{2k-1}x_{2k})$ over the function
field
$\CC (t_{ij})$. The usual linear algebra construction of a corresponding symplectic basis $u_1, \: v_1, \: u_2, \: v_2, \dots , u_m, \: v_m \in
\CC (t_{ij})^{2m}$ goes through: start with any $u_1\neq 0$, find $v_1$ with $\langle u_1, \: v_1\rangle =1$, put $H:=\mathrm{span}(u_1, \: v_1)$,
decompose $\CC (t_{ij})^{2m}= H \oplus H^{\perp }$, and continue with the symplectic form $\langle \cdot , \cdot \rangle_{H^{\perp }}$ in the same
way. This fails for orthogonal groups since one cannot take square roots rationally.

\

(d): The assertion follows immediately by the application of the techniques in part (b). If $P\to X$ is a principal $G$-bundle, then the
associated $G/H$-principal bundle $P\times^G (G/H) \to X$ has locally around any point of $X$ a section because $G/H$ is special. $P$ is an
$H$-principal bundle over
$P\times^G (G/H)$ and pulling back via the section one obtains an $H$-principal bundle $Q$ locally around any point of $X$ which is Zariski
locally trivial because
$H$ is special; and
$P$ is just
$Q\times^H G$, thus is Zariski locally trivial, too.\\
The second assertion follows because a connected linear algebraic group is an extension of its reductive part by the unipotent radical (connected
solvable), and the reductive part an extension of the semi-simple part by a torus.
\end{proof}

\begin{remark}\xlabel{rGrothspecial}
Grothendieck \cite{Groth58} has shown that the only special semi-simple groups are the products of the groups of type $\mathrm{SL}_n (\CC )$ and
$\mathrm{Sp}_n (\CC )$. Serre \cite{Se58} has shown that any special algebraic group is linear and connected.
\end{remark}

If $X$ is a $G$-variety, $G$ a linear algebraic group, one needs a practically verifiable condition when $X\to X/G$ is generically a $G$-principal bundle.

\begin{definition}\xlabel{dfreeaction}
\begin{itemize}
\item[(1)]
The action of $G$ on $X$ is called \emph{free} if the morphism $G\times X \to X \times X$, $(g,\: x)\mapsto (gx, \: x)$ is a closed embedding.
\item[(2)]
$G$ is said to act on $X$ with \emph{trivial stabilizers} if for each point $x\in X$ the stabilizer $G_x$ of $x$ in $G$ is reduced to the identity.
\end{itemize}
\end{definition}

Unfortunately, (1) and (2) are not equivalent. Mumford (\cite{Mum}, Ex. 0.4) gives an example of an action of the group $\mathrm{SL}_2 (\CC )$ on
a quasi-projective variety with trivial stabilizers, but which is not free. However, when for each
$x\in U$,
$U\subset X$ some open dense set, the stabilizers are trivial, we will nevertheless sometimes say that $G$ acts \emph{generically freely} since
this has become standard terminology. We also say more accurately that $G$ acts with \emph{generically trivial stabilizers}.\\
Despite the presence of the subtlety which is displayed in Mumford's example, one has the following result. 

\begin{theorem}\xlabel{tAndrew}
Let $G$ be a connected linear algebraic group acting on a variety $X$ with trivial stabilizers and let $X \to X/G$ be a geometric quotient. Then
there is an open dense $G$-invariant subset $U\subset X$ such that $U \to U/G$ is a $G$-principal bundle in the \'{e}tale topology.
\end{theorem}
\begin{proof}
We use a Seshadri cover (\cite{Sesh72}, \cite{BB} \S 8.4): given a connected linear algebraic group which acts on a variety $X$ with finite
stabilizers, there exists a finite morphism $\kappa \, :\, X_1 \to X$ with the following properties:
\begin{itemize}
\item
$X_1$ is a normal variety and $\kappa$ a (ramified) Galois cover with Galois group $\Gamma$ acting on $X_1$.
\item
There exists a free action of $G$ on $X_1$, commuting with the action of $\Gamma$, such that $\kappa$ is $G$-equivariant.
\item
There exists a good geometric quotient $\pi\, :\, X_1 \to X_1/G$ with $X_1/G$ a prevariety (not necessarily separated), and $\pi$ is a Zariski
locally trivial $G$-principal bundle.
\end{itemize}
"Good" geometric quotient means that $\pi$ is affine. Now assume that in our original situation $V\to V/G$ is a geometric quotient, $V\subset X$
open. Shrinking $V/G$ (and $V$) we can find a $G$-invariant open set $U\subset X$ such that in the diagram
\[
\begin{CD}
\kappa^{-1}(U) @>{\kappa}>>  U\\
@V{\pi }VV                        @VVV\\
\kappa^{-1}(U)/G @>{\bar{\kappa}}>>        U/G
\end{CD}
\]
all arrows are geometric quotients, $\kappa^{-1}(U)/G$ is a variety, and $\bar{\kappa}$ is \emph{\'{e}tale}. It follows that $U \to
U/G$ is a $G$-principal bundle in the \'{e}tale topology since it becomes one (even a Zariski locally trivial one) after the \'{e}tale base change
$\bar{\kappa }$ (the hypothesis that $G$ acts with trivial stabilizers on $X$ implies that the above diagram is a fibre product).
\end{proof}

\begin{corollary}\xlabel{cstabratspecialgroups}
Let a connected linear algebraic group $G$ whose semi-simple part is a direct product of groups of type $\mathrm{SL}$ or $\mathrm{Sp}$ act on a
rational variety $X$ with generically trivial stabilizers. Then $X/G$ is stably rational.
\end{corollary}

\begin{proof}
By Rosenlicht's Theorem \ref{Rosenlicht} and the preceding Theorem \ref{tAndrew}, there is a nonempty open $G$-invariant subset $U\subset X$ such
that a geometric quotient $U \to U/G$ exists and is a $G$-principal bundle in the \'{e}tale topology. This is Zariski locally trivial by Theorem
\ref{tspecialgroups}, (d). Thus to conclude the proof it suffices to remark that a connected linear algebraic group (over $\CC $) is a rational
variety: take a Borel subgroup $B$ in $G$ and consider the $B$-principal bundle $G \to G/B$ over the (rational) flag variety $G/B$. Note that $B$
is rational since it is a successive extension of groups $\mathbb{G}_a$ and $\mathbb{G}_m$.
\end{proof}

However, Corollary \ref{cstabratspecialgroups}, as it stands, is not applicable when we consider for example the action of $\mathrm{PGL}_3 (\CC )$
on the space $\PP (\mathrm{Sym}^4 (\CC^3 )^{\vee })$ of plane quartics. One has the following easy extension.

\begin{corollary}\xlabel{cspecialandtorus}
Let $V$ be a linear representation of a connected linear algebraic group with semi-simple part a direct product of groups $\mathrm{SL}$ and
$\mathrm{Sp}$. Suppose that the generic stabilizer of $G$ in $V$ is trivial. Then $\PP (V) \dasharrow \PP (V)/G$ has a rational section.
\end{corollary}

\begin{proof}
By Rosenlicht's Theorem \ref{Rosenlichtsections}, we see that $V/G \dasharrow \PP (V) /G$ has a rational section, and composing with a rational
section of $V \dasharrow V/G$ and the projection $V \dasharrow \PP(V)$, we obtain a rational section of $\PP (V) \dasharrow \PP (V)/G$.
\end{proof}

\section{Cones and homogeneous bundles}

Over the last decades a variety of different techniques have been developed to make progress on the rationality problem \ref{problemrat} in
certain special cases. These methods will be discussed in the next chapter. However, one may be left with the impression that this is a somewhat
incoherent arsenal of tricks, and no conceptual framework has yet been found which Problem \ref{problemrat} fits into. The purpose of this section
is therefore to discuss certain concepts which seem to have an overall relevance to Problem \ref{problemrat}, and in particular, show how the
Hesselink stratification of the nullcone and the desingularizations of the strata in terms of homogeneous vector bundles give a strategy for
approaching Problem \ref{problemrat}.

\subsection{Torus orbits and the nullcone}

$G$ will be a reductive linear algebraic group throughout this section, $T\subset G$ a fixed maximal torus. $X^{\ast }(T)$ is the group of
characters $\chi \, :\, T \to \CC^{\ast}$ of $T$, and $X_{\ast }( T)$ the group of cocharacters or one-parameter subgroups (1-psg) $\lambda\,
:\, \mathbb{C}^{\ast} \to T$ of $T$. There is the perfect pairing of lattices
\[
\langle \: , \: \rangle \, :\, X^{\ast }(T) \times X_{\ast}(T) \to \ZZ
\]
with $\langle \chi , \lambda \rangle$ defined by $\chi (\lambda (s))= s^{\langle \chi , \lambda \rangle }$ for $s\in \CC^{\ast }$. Let $V$ be a $T$-module
with weight space decomposition
\[
V = \bigoplus_{\chi \in X^{\ast }(T)} V_{\chi }\, .
\]

\begin{definition}\xlabel{dsupport}
Let $v=\sum_{\chi \in X^{\ast }(T)} v_{\chi }$ be the decomposition of a vector $v\in V$ with respect to the weight spaces of $V$. Let
$\mathrm{supp} (v)$ denote the set of those $\chi$ for which $v_{\chi}\neq 0$ (the \emph{support} of $v$), let $\mathrm{Wt} (v)$ be the convex hull
of
$\mathrm{supp} (v)$ in the vector space $X^{\ast }(T)_{\QQ}:= X^{\ast }(T)\otimes_{\ZZ } \QQ$ (the \emph{weight polytope} of $v$), and $C(v)$ the
closed convex cone generated by the vectors in $\mathrm{supp} (v)$ in $X^{\ast }(T)_{\QQ}$.
\end{definition}

The geometry of $T$-orbit closures in the affine or projective cases in $V$ resp. $\PP (V)$ is completely encoded in $\mathrm{supp} (v)$ and
$C(v)$ resp. $\mathrm{Wt} (v)$ in the following way.

\begin{theorem}\xlabel{tTorbitclosures}
\begin{itemize}
\item[(1)]
If $F$ is a face of the cone $C(v)$, put \[ v_F := \sum_{\chi \in \mathrm{supp}(v) \cap F} v_{\chi } \, . \] Then the map $F \mapsto T\cdot v_F$ is
a bijection of the set of faces of $C(v)$ and the set of $T$-orbits in $\overline{T\cdot v}$. If $F_1$ and $F_2$ are faces of $C(v)$, then
\[
F_1 \subset F_2 \iff  T\cdot v_{F_1} \subset \overline{T\cdot v_{F_2}} \, .
\]
\item[(2)]
If $v\neq 0$ is in $V$, and $X_v := \overline{T\cdot [v]} \subset \PP (V)$ is the torus orbit closure of $[v]$ in $\PP (V)$, then the $T$-orbits
on $X_v$ are in bijection with the faces of the weight polytope $\mathrm{Wt}(v)$: for any point $[w] \in X_v$, $\mathrm{Wt}(w)$ is a face of
$\mathrm{Wt}(v)$. For $[w_1], \: [w_2] \in X_v$ one has $T\cdot [w_1] \subset \overline{T\cdot [w_2]}$ if and only if $\mathrm{Wt}(w_1) \subset
\mathrm{Wt}(w_2)$.
\end{itemize}
\end{theorem}

See \cite{B-S}, Prop. 7, p. 104, and \cite{GKZ}, Chapter 5, Prop. 1.8, for a proof. \\
Since $G$ is reductive, it is well known that $\CC [V]^G$ is finitely generated and closed orbits are separated by $G$-invariants (e.g.
\cite{Muk}, Thm. 4.51, Thm. 5.3), thus for $v\in V$ there is a unique closed orbit in $\overline{G\cdot v}$. Thus we can state the following
Hilbert-Mumford theorem.

\begin{theorem}\xlabel{tHilbertMumford}
Let $G$ be a reductive group, $V$ a (finite dimensional) $G$-representation, and pick $v\in V$. Let $X$ be the unique closed orbit in
$\overline{G\cdot v}$. Then there is a 1-psg
$\lambda
\, :\,
\CC^{\ast}
\to G$ such that $\lim_{t\to 0}\lambda (t) \cdot v$ exists and is in $X$.
\end{theorem}

\begin{proof}
It follows from part (1) of Theorem \ref{tTorbitclosures} that every torus orbit in the torus orbit closure of $v$ can be reached as the limit of
$v$ under a suitable 1-psg (the 1-psg corresponds to an integral linear form in $\mathrm{Hom} (X^{\ast } (T), \ZZ)$ defining the face $F$
corresponding to the torus orbit we want to reach). Thus one just has to prove that there is a
$g\in G$ and a torus
$T\subset G$ with
$\overline{T\cdot gv}\cap X
\neq
\emptyset$. One has the Cartan decomposition $G=K \cdot T \cdot K$ of $G$ where $K$ is a maximal compact subgroup of $G$, and $T$ the
complexification of a maximal torus in $K$.\\
Suppose that $\overline{T\cdot gv}\cap X = \emptyset$ \emph{for all} $g$. Then, since $T$-invariants separate disjoint closed $T$-invariant subsets, one may
find for each $w\in G\cdot v$ a function $I_w \in \CC [V]^T$ which is identically $0$ on $X$ and takes the value $1$ in $w$. The compact set
$K\cdot v$ is -as it is a subset of $G\cdot v$- covered by the open sets $U_w$ in $V$ where $I_w$ does not vanish, for $w$ running through $G\cdot
v$, hence it is covered by finitely many of them. The sum of the absolute values of the $I_w$'s corresponding to this finite covering family is
then a continuous $T$-invariant function $s$ which is strictly positive on the compact $K\cdot v$, but $0$ on $X$. But then $s$ is still strictly
positive on $\overline{TK\cdot v}$ which is a contradiction: since $\overline{G\cdot v}= K\cdot \overline{TK\cdot v}$ as $K$ is compact, the fact
that $\overline{TK\cdot v}$ does not meet $X$ implies that $\overline{G\cdot v}$ does not meet $X$ which is false.
\end{proof}

\begin{definition}\xlabel{dstability}
Suppose $v\in V$ is a vector in the $G$-representation $V$.
\begin{itemize}
\item[(1)]
$v$ is \emph{unstable} if $0\in \overline{G\cdot v}$. The set of these is denoted by $\mathfrak{N}_G(V)$ (the \emph{nullcone}).
\item[(2)]
$v$ is \emph{semistable} if $v$ is not unstable.
\item[(3)]
$v$ is \emph{stable} if $G\cdot  v$ is closed in $V$ and the stabilizer subgroup $G_v\subset G$ of $v$ is finite.
\item[(4)]
$v$ is called $T$\emph{-unstable} if $0\in \overline{T\cdot v}$, and the set of $T$-unstable elements is denoted by $\mathfrak{N}_T (V)$ (the
\emph{canonical cone}).
\end{itemize}
\end{definition}

\begin{theorem}\xlabel{tinstability}
\begin{itemize}
\item[(1)]
The set $\mathfrak{N}_G (V)$ is defined by the vanishing of all invariants in $\CC [V]^G$ of positive degree, $\mathfrak{N}_T(V)$ is defined by
the vanishing of all $T$-invariants in $\CC [V]^T$ of positive degree; this gives $\mathfrak{N}_G (V)$ and $\mathfrak{N}_T(V)$ scheme structures.
\item[(2)]
A vector $v\in V$ is $T$-unstable if and only if $0\notin \mathrm{Wt}(v)$. One has $G\cdot \mathfrak{N}_T(V) = \mathfrak{N}_G (V)$. Hence
$\mathfrak{N}_G (V)$ consists of those vectors $v$ such that the orbit $G\cdot v$ contains an element whose weight polytope does not contain $0$.
\end{itemize}
\end{theorem}

\begin{proof}
(1) follows immediately from the fact that $G$- or $T$-invariants separate closed orbits. The first assertion of (2) is a consequence of Theorem
\ref{tTorbitclosures}, (1). The fact that $G\cdot \mathfrak{N}_T(V) = \mathfrak{N}_G (V)$ follows from Theorem \ref{tHilbertMumford}.
\end{proof}

\begin{remark}\xlabel{rschemestructurenullcone}
In general, $\mathfrak{N}_G(V)$ need neither be irreducible nor reduced nor equidimensional: in \cite{Po92}, Chapter 2 \S 3, it is shown that for
the representation of $\mathrm{SL}_2 (\CC )$ in the space of binary sextics $\mathrm{Sym}^6 (\CC^2 )^{\vee}$, the nullcone is not reduced; for the
representation of $\mathrm{SL}_3 (\CC )$ in the space of ternary quartics $\mathrm{Sym}^4 (\CC^3 )^{\vee }$ the nullcone has two irreducible
components of dimensions $10$ and $11$, respectively (\cite{Hess79}, p. 156).
\end{remark}

\begin{remark}\xlabel{rgraphicalinstab}
Theorem \ref{tinstability} (2) gives a convenient graphical way for the determination of unstable vectors known to Hilbert (\cite{Hil93}, \S 18):
for example, if $\mathrm{SL}_2 (\CC )$ acts on binary forms of degree $d$, $\mathrm{Sym}^d (\CC^2 )^{\vee }$, in variables $x$ and $y$, $T
\subset
\mathrm{SL}_2 (\CC )$ is the standard torus
\[
T = \left( \begin{array}{cc} t & 0 \\ 0 & t^{-1} \end{array} \right), \; t\in \CC^{\ast },
\]
and $\epsilon_1$ is the weight
\[
\epsilon_1 \left( \left( \begin{array}{cc} t & 0 \\ 0 & t^{-1} \end{array} \right) \right) := t ,
\]
then the weight spaces are spanned by $x^{k} y^{d-k}$ which is of weight $(d-2k)\epsilon_1$, for $k=0, \dots , d$. The support of a binary degree
$d$ form thus does not contain the origin if and only if it is divisible by $x^{[d/2]+1}$ or $y^{[d/2]+1}$, so that the unstable binary degree $d$
forms are just those which have a zero of multiplicity $\ge [d/2] +1$.\\
Turning to ternary forms, the representation of $\mathrm{SL}_3 (\CC )$ in $\mathrm{Sym}^d (\CC^3)^{\vee}$ (coordinates $x$, $y$, $z$), one can use
"Hilbert's triangle": consider in the plane an equilateral triangle $ABC$ with barycenter the origin $0$. Let $\epsilon_1, \: \epsilon_2, \: \epsilon_3$ be
the vectors pointing from $0$ to $A, \: B, \: C$. Points of the plane with integer coordinates with respect to the basis $\epsilon_1$, $\epsilon_2$ are thus
identified with the character lattice of the standard maximal torus of
$\mathrm{SL}_3 (\CC )$. For the monomial 
$x^ay^bz^c$,
$a+b+c=d$, the point in the plane $-a\epsilon_1-b\epsilon_2-c\epsilon_3$ then represents the associated weight; pick a line $l$ in this
plane not passing through zero, and let $H_l$ be the corresponding closed half-space in the plane not containing $0$. Then a ternary degree
$d$ form
$f$ is unstable if, after a coordinate change, it may be written as a sum of monomials whose weights lie entirely in $H_l$ for some line $l$. In this way
it is possible to obtain finitely many representatives for all possible types. 
\end{remark}

\subsection{Stratification of the nullcone and rationality}

The importance of the nullcone $\mathfrak{N}_G (V)$ for us derives from the fact that it contains a lot of rational subvarieties which are
birational to homogeneous vector bundles over generalized flag varieties whose fibres are linearly embedded in $V$. We are now going to describe
this.\\
Let $G$ be reductive as before, $T\subset G$ a maximal torus. Let $W:= N_G (T) / Z_G (T)$ be the Weyl group which acts on $T$ by conjugation,
hence on $X^{\ast }(T)$. Choose a $W$-invariant scalar product $\langle \cdot , \cdot \rangle$ on $X^{\ast }(T)_{\QQ }$ which takes integral
values on $X^{\ast }(T)$.\\
Denote by $R \subset X^{\ast }(T)$ the set of roots, the set of nonzero weights of $T$ in the adjoint representation of $G$ on its Lie algebra
$\mathfrak{g}$. For $\alpha \in R$ denote by $U_{\alpha}$ the root subgroup corresponding to $\alpha$, i.e. the unique connected $T$-invariant
unipotent subgroup $U_{\alpha}\subset G$ with Lie algebra the one dimensional root subspace $\mathfrak{g}_{\alpha }$. 

\begin{definition}\xlabel{dStructNullcone}
\begin{itemize}
\item[(1)]
For $c\in X^{\ast} (T)_{\QQ }$ denote by $P_c$ the subgroup of $G$ generated by the torus $T$ and the root groups $U_{\alpha }$ for $\alpha \in R$
with $\langle \alpha , \: c \rangle \ge 0$. Furthermore, $U_c$ resp. $L_c$ will be the subgroups generated by $U_{\alpha}$ with $\langle \alpha ,
\: c \rangle > 0$ resp. by the torus $T$ and the root subgroups $U_{\alpha }$ with $\langle \alpha , \: c\rangle =0$.\\
\item[(2)]
For $c \in X^{\ast } (T)_{\QQ }$, $c\neq 0$, one denotes by $H^+ (c)$ the half space $\{ x \in X^{\ast } (T)_{\QQ } \, | \, \langle x-c,\: c\rangle \ge 0
\}$ in $X^{\ast }(T)_{\QQ }$ and by $H^0 (c) = \{ x \in X^{\ast } (T)_{\QQ } \, | \, \langle x-c, \: c \rangle =0 \}$ its bounding hyperplane.
\item[(3)]
For a subset $\Sigma \subset X^{\ast } (T)_{\QQ }$ and a finite-dimensional $G$-module $V$, $V_{\Sigma}$ denotes the subspace of $V$ consisting of those
$v\in V$ with
$\mathrm{supp}(v) \subset \Sigma$. Obviously, there are only finitely many such subspaces since the number of weight spaces in $V$ is finite.
\end{itemize}
\end{definition}

Remark that $P_c \subset G$ is a parabolic subgroup containing $T$ (because the subset of those $\alpha \in R$ with $\langle \alpha, \: c\rangle
\ge 0$ contains some basis for the root system $R$ and is closed with respect to addition); it follows from this observation that $L_c$ is a
reductive Levi subgroup of $P_c$ containing $T$, and $U_c$ is the unipotent radical of $P_c$, $P_c = L_c \ltimes U_c$. Moreover:

\begin{lemma}\xlabel{lCoarseStructNullcone}
\item[(1)]
For $c \in X^{\ast }(T)_{\QQ }$, $c \neq 0$, the subspace $V_{H^+ (c)}$ is stable under $P_c$, the subspace $V_{H^0 (c)}$ is stable under the Levi
subgroup $L_c$.
\item[(2)]
$G\cdot V_{H^+ (c)}$ is closed in $V$, and the image of the homogeneous vector bundle $G \times_{P_c} V_{H^+ (c)} \to G/ P_c$ under the
natural $G$-map to $V$. Moreover, the nullcone can be expressed as a union (which is actually finite) of such images:
\[
\mathfrak{N}_G (V) = \bigcup\limits_{c \neq 0,\, c \in X^{\ast�}(T)_{\QQ }} G \cdot V_{H^+ (c)} \, .
\]
\end{lemma}

\begin{proof}
(1) follows from the following well-known fact from the representation theory of reductive groups: if $\chi \in X^{\ast }(T)$,
$v\in V_{\chi�}$,
$\alpha\in R$ and
$g\in U_{\alpha }$, then
\[
g\cdot v -v \in \bigoplus\limits_{l\ge 1} V_{\chi +l\alpha }\, .
\]
It is seen as follows: let $x_{\alpha} \, :\, \mathbb{G}_a \to G$ be the root homomorphism which is an isomorphism onto its image $U_{\alpha }$
and $t\cdot x_{\alpha }(k)\cdot t^{-1} = x_{\alpha }(\alpha (t) k)$ $\forall\, t\in T, \: \forall \, k \in\mathbb{G}_a$. Then $x_{\alpha }(k)\cdot
v$ is a polynomial in $k$ with coefficients in $V$:
\[
x_{\alpha}(k)\cdot v = \sum\limits_{l=0}^N v_l k^l, \; \forall \, k\in \mathbb{C}
\]
and on the one hand, for $t\in T$
\[
t\cdot (x_{\alpha }(k) \cdot v) = \sum\limits_{l=0}^N k^l t\cdot v_l \, ,
\]
whereas on the other hand
\begin{gather*}
t\cdot (x_{\alpha }(k) \cdot v) = (t\cdot x_{\alpha}(k) t^{-1}) (tv)\\
= x_{\alpha}(\alpha (t) \cdot k) \chi (t) v = \sum\limits_{l=0}^N \alpha (t)^l k^l \chi (t) \cdot v_l \, .
\end{gather*}
Equating coefficients in the two polynomials in $k$ yields $\chi (t) \alpha (t)^l\cdot v_l = t \cdot v_l$, all $t\in T$, which together with the
fact that $v=v_0$ (put $k=0$) gives the desired result.

\ 

(2): $G\cdot V_{H^+ (c)}$ is the image of the natural $G$-map $G\times_{P_c} V_{H^+ (c)} \to V$ which factors into the closed
embedding 
$i\, :\, G\times_{P_c} V_{H^+ (c)} \to G/P_c \times V$ given by $i ([(g,\: v)]) := (gP_c, \: g\cdot v)$, followed by the projection $G/P_c \times V
\to V$ onto the second factor which is proper since $G/P_c$ is compact. Thus $G\cdot V_{H^+ (c)}$ is closed in $V$. The last assertion about the
nullcone is an immediate consequence of Theorem \ref{tinstability}, (2).
\end{proof}

We need a criterion for when the natural $G$-map $G\times_{P_c} V_{H^+ (c)} \to V$ is birational onto its image $G\cdot V_{H^+ (c)}$.

\begin{theorem}\xlabel{tstratanullcone}
Let $v \in V$ be a $T$-unstable element. Then the norm $|| \cdot || := \sqrt{ \langle \cdot , \cdot \rangle }$ induced by the $W$-invariant scalar
product on $X^{\ast }(T)_{\QQ }$ achieves its minimum in exactly one point $c$ of $\mathrm{Wt}(v)$. Thus there exists a smallest positive integer
$n$ such that $n\cdot c$ is in $X^{\ast}(T)$. Since $n\cdot c$ is orthogonal to the roots of $L_c$, it extends to a character of $P_c$ and
$L_c$. Let
$Z_c:= (\mathrm{ker}(n\cdot c|_{L_c}) )^{\circ } \subset L_c$ be the corresponding reductive subgroup of $L_c$ with maximal torus $T' =
(\mathrm{ker}(n\cdot c) \cap T)^{\circ}$.\\ The space $V_{H^+ (c)}$ decomposes as $V_{H^+(c)} = V_{H^0 (c)} \oplus \bigoplus_{\chi \in X^{\ast
}(T), \: \chi \in H^+ (c)\backslash H^0(c)} V_{\chi}$. Let $v_0$ be the component of $v$ in $V_{H^0(c)}$ with respect to this decomposition.\\
Suppose that $v_0$ is not in the nullcone of $Z_c$ in $V_{H^0(c)}$. Then the $G$-map
\begin{gather*}
G\times_{P_c} V_{H^+(c)} \to G\cdot V_{H^+(c)} \subset V\\
([(g, v')]) \mapsto g\cdot v'
\end{gather*}
is birational onto $G\cdot V_{H^+(c)}$ and the fibres of the bundle $G\times_{P_c} V_{H^+(c)}$ are embedded as linear subspaces of $V$.
\end{theorem}

\begin{proof}
All faces of the polytope $\mathrm{Wt}(v)$ are described by rational linear equalities and inequalities, and $|| \cdot ||^2$ is
a rational quadratic form on $X^{\ast}(T)_{\QQ }$. By the differential criterion for extrema with boundary conditions, one gets a system of linear
equations defined over $\QQ$ which determine $c$.\\
We now use an argument in \cite{Po-Vi} \S 5, due to Kirwan and Ness. Remark that the set of vectors $u$ in $V_{H^+(c)}$ for which $c$ is the point
of
$\mathrm{Wt}(u)$ closest to the origin and for which the projection $u_0$ onto $V_{H^0(c)}$ is not unstable for $Z_c$ is open in $V_{H^+(c)}$ and
not empty by assumption. We denote this set by $\Omega^+(c)$. Let $g_1, \: g_2 \in G$, $u_1, \: u_2 \in \Omega^+(c)$ be such that $g_1 u_1 = g_2
u_2$. So that $([(g_1, \: u_1)])$ and $([(g_2, \: u_2)])$ in the bundle $G\times_{P_c} V_{H^+(c)}$ map to the same image point in $V$. It thus
suffices to show that for each $g\in G$, $gu_1 =u_2$ implies $g\in P_c$.\\
Look at the Bruhat decomposition $G=U_{P_c} W P_c$ where $U_{P_c}$ is the unipotent radical of $P_c$. We may thus write $g=p_1 w p_2$ with
$p_1, \: p_2\in P_c$, $w$ some (representative of an) element in $W$. We now make the following
\begin{quote}
Claim: For each $p\in P_c$ and every $u\in\Omega^+(c)$, the weight polytope $\mathrm{Wt}(pu)$ contains the point $c$ (one cannot "move
$\mathrm{Wt}(u)$ away from $c$" with elements in $P_c$).
\end{quote}
Assuming this claim for the moment, one can finish the proof as follows: rewriting $g u_1 =u_2$ as $w (p_2 u_1) = p_1^{-1}u_2$ we see that
$w\cdot c \in w \cdot \mathrm{Wt}(p_2u_1) = \mathrm{Wt}(w p_2 u_1) \subset H^+(c)$. But since $c$ is the only element of norm $|| c||$ in
$H^{+}(c)$ and $W$ acts by isometries, we must have $w\cdot c = c$. Thus $P_c = P_{w\cdot c} = w P_c w^{-1}$, whence $w\in N_G(P_c) = P_c$.\\
We turn to the proof of the claim. In fact, the claim is just a reformulation of the property of the $u\in\Omega^+(c)$ to have $Z_c$-semistable
projection $u_0$ onto $V_{H^0(c)}$ and the Hilbert-Mumford criterion in the form of Theorem \ref{tinstability}, (2).  
Recall from above that if $\chi \in X^{\ast }(T)$,
$v\in V_{\chi�}$,
$\alpha\in R$ and
$g\in U_{\alpha }$, then
$g\cdot v -v \in \bigoplus\limits_{l\ge 1} V_{\chi +l\alpha }$, so if $\langle \alpha , c \rangle >0$, then replacing $u$ by $f\cdot u$
for $f \in U_{\alpha}$ gives 
$\mathrm{Wt}(u) \cap H^0 (c)= \mathrm{Wt}(f\cdot u)\cap H^0 (c)$. Thus the weight polytope of $u$ can be moved away from $c$ by an element of
$P_c$ if and only if it can be moved away from $c$ by an element of $L_c$ and $L_c$ is in turn generated by $Z_c$ and a one-dimensional central
subtorus in $L_c$ (central since all roots of $L_c$ are trivial on it). But under the restriction of characters of $T$ to the subtorus $T'$
which is the maximal torus in $Z_c$, $H^0(c)$ maps bijectively onto $X^{\ast}(T')_{\QQ }$ and the point $c$ gets identified with the origin in
$X^{\ast }(T')_{\QQ }$. Thus the weight polytope of $u$ in $X^{\ast}(T)_{\QQ }$ can be moved away from $c$ by $P_c$ if and only if the weight
polytope of the projection $u_0$ of $u$ onto $V_{H^0(c)}$ can be moved away from $0$ in $X^{\ast }(T')_{\QQ }$ by the action of $Z_c$. Our
assumption that $u_0$ be $Z_c$-semistable exactly prevents this possibility.
\end{proof}

\begin{remark}\xlabel{rStructureStrata}
Note that, though the way it was stated, Theorem \ref{tstratanullcone} involves the choice of a $T$-unstable element $v$, and $c\in
X^{\ast}(T)_{\QQ }$ is afterwards determined as the element of $\mathrm{Wt}(v)$ of minimal length, the \emph{only important requirement} is that
$\mathfrak{N}_{Z_c} (V_{H^{0}(c)}) \neq V_{H^0(c)}$. In fact, if $c\in X^{\ast}(T)_{\QQ }$, $c\neq 0$, is any element with this property, then
the set $\Omega^+(c) \subset V_{H^+(c)}$ of vectors whose projection unto $V_{H^0(c)}$ is $Z_c$-semistable is nonempty, and the weight polytopes of all these
vectors must automatically contain $c$ then, and are contained in $H^+(c)$, hence $c$ is the vector of minimal distance to the origin in all those
weight polytopes. 
\end{remark}

\begin{definition}\xlabel{dstratifyingelement}
An element $c\in X^{\ast}(T)_{\QQ}$, $c\neq 0$, with $\mathfrak{N}_{Z_c} (V_{H^{0}(c)}) \neq V_{H^0(c)}$ is called \emph{stratifying}.
\end{definition}

The finer structure of the nullcone is described in

\begin{theorem}\xlabel{tpartitionnullcone}
Let $c\in X^{\ast }(T)_{\QQ }$ be a stratifying element. Then $G\times_{P_c} \Omega^+(c) \to G\cdot \Omega^+(c)=:\mathfrak{S}(c)$ is an
isomorphism, and $G\cdot V_{H^+(c)}\subset \mathfrak{N}_G(V)$ is the closure of $\mathfrak{S}(c)$. We call $\mathfrak{S}(c)$ a nonzero stratum of
$\mathfrak{N}_G (V)$. One has $\mathfrak{S}(c_1) = \mathfrak{S}(c_2)$ if and only if $c_1$ and $c_2$ are in the same $W$-orbit, and
$\mathfrak{N}_G (V)$ is a finite disjoint union of $\{ 0\}$ and the nonzero strata.
\end{theorem}

See \cite{Po-Vi}, \S 5, \cite{Po03}. Through the paper \cite{Po03}, an algorithm -using only the configuration of weights with multiplicities of
$V$ and the roots of $G$ in $X^{\ast }(T)_{\QQ }$- is now available to determine completely the family of stratifying elements resp. strata. We
pass over all this, since for our immediate purposes, Theorem \ref{tstratanullcone} is sufficient.\\
We will now describe how Theorem \ref{tstratanullcone} can be used to develop a general technique for approaching the rationality problem
\ref{problemrat}. This was suggested in \cite{Shep89} and we will develop it in greater detail.\\

\begin{theorem}\xlabel{tMainTrick}
Let $\Gamma \subset G$ be connected reductive groups, $V$ a $G$-module, and $M$ a $\Gamma$-submodule of $V$ which is contained in the nullcone
$\mathfrak{N}_G (V)$ of $G$ in $V$. Let $\mathfrak{S}(c)$, $0\neq c\in X^{\ast }(T)_{\QQ }$, $T\subset G$ a maximal torus, be a stratum of
$\mathfrak{N}_G (V)$. Let
\[
G\times_{P_c} V_{H^+(c)}  \to G \cdot V_{H^+(c)} = \overline{\mathfrak{S}(c)} \subset V
\]
be the associated desingularization of $\overline{\mathfrak{S}(c)}$ by the homogeneous vector bundle
\[
\begin{CD}
G\times_{P_c} V_{H^+(c)} @> \pi >> G/P_c \, .
\end{CD}
\]
Assume:
\begin{itemize}
\item[(a)]
$\mathfrak{S}(c) \cap M$ is dense in $M$ and the rational map $\pi\, :\, \PP (M) \dasharrow G/P_c$ induced by the bundle projection is dominant.
\item[(b)]
$(G/P_c)/ \Gamma$ is stably rational in the sense that $(G/P_c)/ \Gamma \times \PP^r$ is rational for some $r\le \dim \PP (M) - \dim G/P_c$.
\item[(c)]
Let $Z$ be the kernel of the action of $\Gamma$ on $G/P_c$: assume $\Gamma /Z$ acts generically freely on $G/P_c$, $Z$ acts trivially on $\PP (M)$, and there
exists a
$\Gamma /Z$-linearized line bundle
$\mathcal{L}$ on the product $\PP (M) \times G/P_c$ cutting out $\mathcal{O}(1)$ on the fibres of the projection to $G/P_c$.
\end{itemize}
Then $\PP (M) / \Gamma$ is rational.
\end{theorem}

\begin{proof}
Let $Y:= G/P_c$, $X:=$the (closure of) the graph of $\pi$, $p\, :\, X\to Y$ the restriction of the projection which (maybe after shrinking $Y$) we
may assume (by (a)) to be a projective space bundle for which $\mathcal{L}$ is a relatively ample bundle cutting out $\mathcal{O}(1)$ on the
fibres. The main technical point is the following result from descent theory (\cite{Mum},
\S 7.1): by Theorem \ref{tAndrew} there are nonempty open subsets $X_0\subset X$ and $Y_0 \subset Y$ such that we have a fibre product square with the
bottom horizontal arrow a $\Gamma /Z$-principal bundle:
\[
\begin{CD}
X_0 @>>> X_0 / (\Gamma /Z)\\
@V{p}VV        @V{\bar{p}}VV\\
Y_0 @>>> Y_0/(\Gamma /Z)
\end{CD}
\]
and by \cite{Mum}, loc. cit., $\mathcal{L}$ descends to a line bundle $\bar{\mathcal{L}}$ on $X_0/(\Gamma /Z)$ cutting out $\mathcal{O}(1)$ on the fibres of
$\bar{p}$. Hence $\bar{p}$ is also a Zariski locally trivial projective bundle (of the same rank as $p$). 
By (b), it now follows that $\PP (M)/\Gamma$ is rational.
\end{proof}

\begin{example}\xlabel{eTwoForms}
Let $E$ be a complex vector space of odd dimension $n$ ($n\ge 3$) and consider the action of the group $G= \mathrm{SL} (E)$ on $V = \Lambda^2 (E)$.
We choose a basis $e_1, \dots , e_n$ of $E$ so that $\mathrm{SL} (E)$ is identified with the group $\mathrm{SL}_n (\CC )$ of $n\times n$ matrices
of determinant one. Let 
\[
T = \mathrm{diag} (t_1, \dots , t_n) , \quad t_i \in \CC , \quad \prod\limits_{i=1}^n t_i =1
\]
be the standard diagonal torus, and denote by $\epsilon_i \in X^{\ast}(T)$ the $i$th coordinate function of $T$
\[
\epsilon_i (\mathrm{diag} (t_1, \dots , t_n)) = t_i, \quad i=1, \dots , n \, .
\]
In $\RR^n$ with its standard scalar product consider the hyperplane $H := \langle (1,1, \dots , 1) \rangle^{\perp}$. We make the identifications
\begin{gather*}
\epsilon_1 = (n-1, \: -1, \: -1, \dots , -1), \quad \epsilon_2 = (-1, \: n-1, \: -1, \dots , -1) , \dots ,\\ \epsilon_n = (-1, \: -1, \: -1, \dots
, n-1) \quad \mathrm{(then }\; \epsilon_i \in H \quad \forall\, i=1, \dots , n \mathrm{)}
\end{gather*}
whence $X^{\ast}(T)\otimes\RR$ becomes identified with $H$, $X^{\ast} (T)$ being the subset of vectors $a_1 \epsilon_1 + \dots + a_n \epsilon_n$
with $a_i \in \ZZ$ for all $i$, and $X^{\ast}(T)_{\QQ }$ is similarly defined by the condition $a_i \in\QQ$. The restriction of the standard
Euclidean scalar product on $\RR^n$ to $H$ is then a $W$-invariant scalar product, integral on $X^{\ast}(T)$. Denote it by $\langle \cdot , \cdot
\rangle$. Note that the Weyl group $W$ is the symmetric group $\mathfrak{S}_n$ on $n$ letters and acts by permuting the $\epsilon_i$. The roots of
$(G,\: T)$ are then
\[
\alpha_{ij} := \epsilon_i - \epsilon_j , \quad  1\le i, \: j \le n , \; i \neq j 
\]
with corresponding root subgroups
\[
U_{\alpha_{ij}} = \left\{ A \in \mathrm{SL}_n (\CC) \, |\, A = \mathrm{Id} + r \cdot E_{ij} \right\}, 
\]
with $E_{ij}$ the $n\times n$ elementary matrix with a single nonzero entry, namely $1$, in the $(i,\: j)$-spot. Thus the $\epsilon_i$ form the
vertices of a simplex in $H\simeq X^{\ast }(T) \otimes\RR $ and the roots are the pairwise differences of the vectors leading from the
origin to the vertices. The weights of $T$ in $\Lambda^2 (E)$ are obviously
\[
\pi_{kl} = \epsilon_k + \epsilon_l , \quad 1 \le k < l \le n, \quad  V_{\pi_{kl}} = \CC\cdot ( e_k \wedge e_l )\, .
\]
Define an element $c \in X^{\ast}(T)_{\QQ }$ by
\[
c := \frac{2}{n-1} (1, \: 1, \: 1, \dots , \: 1, \: -(n-1)) \in H \, .
\]
We consider as above the affine hyperplane $H^0(c)$ perpendicular to $c$ and passing through $c$, and the positive half space $H^+(c)$ it defines.
The following facts concerning the relative position of $H^0(c)$ and the weights $\pi_{kl}$ and roots $\alpha_{ij}$ are easily established by
direct calculation:
\begin{align*}
\bullet &\;\pi_{kl} \in H^+ (c ) \iff \pi_{kl} \in H^0 (c ) \iff 1\le k < l \le n-1 , \\
\bullet &\; \langle \alpha_{ij} , \: c \rangle =0 \iff  1 \le i, \: j \le n-1, \: i\neq j \, ,\\
\bullet &\; \langle \alpha_{ij} , \: c \rangle > 0 \iff 1\le i < j = n \, .
\end{align*}
Note that then in the above notation one has for the group $P_c$
\begin{gather*}
P_c = \left\{  \left( \begin{array}{cc} M & a\\ 0 & b  \end{array} \right) \in \mathrm{SL}_n (\CC ) \, : \,  M \in \CC^{(n-1)\times (n-1)}, \: a\in
\CC^{(n-1)\times 1} , \right. \\ \left.  \: b \in \CC , \: 0\in \CC^{1\times (n-1)}
\right\}\, .
\end{gather*}
Similarly, the reductive group $L_c$ is
\[
L_c = \left\{  \left( \begin{array}{cc} M & 0^t\\ 0 & b  \end{array} \right) \in \mathrm{SL}_n (\CC ) \, : \,  M \in \mathrm{GL}_{n-1} (\CC ), \:
 b \in \CC^{\ast } , \: 0\in \CC^{1\times (n-1)}
\right\}\, .
\]
Now $\frac{n-1}{2} c = \epsilon_1 + \dots + \epsilon_{n-1}$. This extends to the character of $L_c$ which maps an element of $L_c$ to the
determinant of $M$. Hence the group $Z_c$ is
\[
Z_c = \left\{  \left( \begin{array}{cc} M & 0^t \\ 0 & 1  \end{array} \right) \in \mathrm{SL}_n (\CC ) \, : \,  M \in \mathrm{SL}_{n-1}(\CC ), \:
 0\in \CC^{1\times (n-1)} \right\} \, .
\]
The action of $Z_c$ on $V_{H^0 (c)}$ is equivalent to the standard action of $\mathrm{SL}_{n-1}(\CC )$ on $\Lambda^2 (\CC^{n-1})$ and the nullcone
for this action is not the whole space (there exists the Pfaffian). Hence $c$ is stratifying.\\
The flag variety $G/P_c$ can be identified with the Grassmannian $\mathrm{Grass} (n-1, E)$ of $n-1$-dimensional subspaces in $E$, or dually, $\PP
(E^{\vee })$, the projective space of lines in $E^{\vee }$. The open set $G\cdot \Omega^+(c)$ in $\Lambda^2 (E)$ can be identified with the two
forms  represented by skew-symmetric matrices of maximal rank $n-1$. Every vector in $\Lambda^2 (E)$ is unstable. If we view $\Lambda^2(E)$ as
$\Lambda^2 (E^{\vee })^{\vee }$, i.e. skew-forms on $E^{\vee }$, then the bundle projection
\[
G\times_{P_c} \Omega^{+}(c) \to G/P_c
\]
is identified with the map which assigns to a skew-form $\omega$ of maximal rank its image under the linear map
\begin{gather*}
E^{\vee } \to E , \quad  e \mapsto \omega (e, \cdot ) \, ,
\end{gather*}
(an element of $\mathrm{Grass}(n-1, E)$), or dually, its kernel in $\PP (E^{\vee} )$. The associated method for proving rationality is called
\emph{the 2-form trick} and appears in \cite{Shep}, Prop.8.
\end{example}

\begin{example}\xlabel{eDoubleBundle}
Let $E$ and $F$ be complex vector spaces with $\dim E =: n > \dim F =: m$. As in the previous example, for the action of $G =
\mathrm{SL} (E) \times \mathrm{SL} (F)$ on $V= \mathrm{Hom} (E ,\: F)$, every vector is unstable (there is a dense orbit). Choose bases $e_1,
\dots ,\:  e_n$ for $E$ and $f_1, \dots , \: f_m$ for $F$, so that $\mathrm{SL} (E) \times \mathrm{SL} (F) \simeq \mathrm{SL}_{n} (\CC )
\times
\mathrm{SL}_m (\CC )$. If $T^E$ resp. $T^F$ denote the standard maximal tori of diagonal matrices in $\mathrm{SL}_{n} (\CC )$ resp.
$\mathrm{SL}_m (\CC )$, then
\[
T = T^E \times T^F
\]
is a maximal torus of $G$. For $(\mathrm{SL}_{n} (\CC ), \: T^E)$ and $(\mathrm{SL}_m (\CC ), \: T^F)$, we use the definitions and concrete
realization of weight lattices as in Example \ref{eTwoForms}, except that we endow now all objects with superscripts $E$ and $F$  to indicate
which group we refer to: thus we write, for example, $\epsilon_i^E$, $\epsilon_j^F$, $W^E$, $X^{\ast}(T^E)\otimes\RR \simeq H^E$, and so forth.\\
Then we have $X^{\ast}(T) = X^{\ast} (T^E ) \times X^{\ast}(T^F)$ and we may realize $X^{\ast}(T)\otimes \RR$ as 
\[
H := H^E \times H^F \subset \RR^{n} \times \RR^m \simeq \RR^{n+m}
\]
with scalar product
\[
\langle (x_1, \: y_1) , \: (x_2, \: y_2) \rangle := \langle x_1 , \: x_2 \rangle^E + \langle y_1, \: y_2 \rangle^F, \quad x_{1}, \: x_2 \in H^E,
\: y_1, \: y_2 \in H^F \, ,
\]
which is invariant under the Weyl group $W= W^E \times W^F\simeq \mathfrak{S}_{n}\times \mathfrak{S}_m$ of $G$ acting by permutations of the
$\epsilon_i^E$ and $\epsilon_j^F$ separately. The weights of $V$ with respect to $T$ are given by
\[
\pi_{kl} := ( -\epsilon_k^E, \: \epsilon_l^F), \: 1\le k \le n , \: 1 \le l\le m, \quad V_{\pi_{kl}}= \CC \cdot (e_k^{\vee} \otimes f_l )
\]
where under the isomorphism $E^{\vee }\otimes F \simeq \mathrm{Hom} (E, \: F)$, the vector $e_k^{\vee}\otimes f_l$ corresponds to a matrix with
only one nonzero entry $1$ in the $(l, \: k)$ position. Note that the $\pi_{kl}$ form the vertices of a polytope in $H$ which is the product of
two simplices. The roots of $(G,\: T)$ in $H$ are the vectors
\[
(\alpha_{pq}^E, \: 0), \: 1\le p,\: q\le {n}, \: p\neq q, \quad (0, \: \alpha_{rs}^F), \: 1\le r,\: s\le m, \: r\neq s
\]
(the disjoint union of the root systems of $\mathrm{SL}_{n}(\CC )$ and $\mathrm{SL}_m (\CC )$ in the orthogonal subspaces $H^E$ and $H^F$), and
the root subgroups in $G$ are then
\[
U_{(\alpha_{pq}^E, \: 0)} = U_{\alpha_{pq}^E} \times \{ \mathrm{Id}_m \}, \quad U_{(0, \: \alpha_{rs}^F)} = \{ \mathrm{Id}_{n} \} \times
U_{\alpha_{rs}^F} \, .
\]
Now define $c\in X^{\ast}(T)_{\QQ }$ by
\[
c:= \left( -\frac{1}{m}(n-m, \: n-m, \dots , \: n-m, \: -m, \: -m, \dots , \: -m), \; 0 \right) \in H^E \times H^F
\]
where there are $m$ entries with value $n-m$ followed by another $n-m$ entries with value $-m$ in the row vector in the first component. In fact,
$c = -\frac{1}{m} ((\epsilon_1^E, \: 0)  + \dots + (\epsilon_m^E, \: 0) )$. The following facts are easily verified by direct computation:
\begin{align*}
\bullet &\pi_{kl} \in H^{0}(c) \iff \pi_{kl} \in H^+ (c) \iff 1 \le k \le m \, ,\\
\bullet &\langle (0, \alpha^F_{rs} ) , \: c \rangle = 0 \: \forall\,  r,\: s, \: 1\le r, \: s \le m, r\neq s \, , \\
\bullet &\langle (\alpha^E_{pq}, 0) , \: c \rangle = 0 \iff 1\le p,\: q \le m, \: p\neq q \; \mathrm{or}\; m+1 \le p,\: q \le n, \: p\neq q \, ,\\
\bullet &\langle (\alpha^E_{pq}, 0) , \: c \rangle > 0 \iff 1\le q \le m, \: m+1 \le p \le n\, .
\end{align*}
Hence we get
\begin{gather*}
P_c = \left\{ \left( \begin{array}{cc} A & 0\\ B & C     \end{array}\right) \in \mathrm{SL}_n (\CC ) \, : \, A\in \mathbb{C}^{m\times m}, 
 \: B
\in
\CC^{(n-m)\times m}, \: C\in \CC^{(n-m)\times (n-m)} 
\right\}\\
\times
\mathrm{SL}_m (\CC ) \, ,
\end{gather*}
and
\begin{gather*}
L_c = \left\{ \left( \begin{array}{cc} A & 0\\ 0 & C     \end{array}\right) \in \mathrm{SL}_n (\CC ) \, : \, A\in \mathrm{GL}_m (\CC ), 
 \: C\in \mathrm{GL}_{n-m} (\CC ) 
\right\} \times \mathrm{SL}_m (\CC )\, , \\
Z_c = \left\{ \left( \begin{array}{cc} A & 0\\ 0 & C     \end{array}\right) \in \mathrm{SL}_n (\CC ) \, : \, A\in \mathrm{SL}_m (\CC ), 
 \: C\in \mathrm{SL}_{n-m} (\CC ) 
\right\} \times \mathrm{SL}_m (\CC )\, ,
\end{gather*}
since $- m\cdot c$ extends to the character of $L_c$ which maps an element of $L_c$ to the determinant of $A$. The action of $Z_c$ on $V_{H^0(c)}$
is equivalent to the standard action of  $\mathrm{SL}_m (\CC ) \times \mathrm{SL}_m (\CC )$ on $\mathrm{Hom} (\CC^m, \: \CC^m)$, whence
$\mathfrak{N}_{Z_c} (V_{H^0(c)}) \neq V_{H^0(c)}$ since an endomorphism has a determinant. Thus $c$ is stratifying.\\
The flag variety $G/P_c$ is isomorphic to the Grassmannian $\mathrm{Grass}(n-m, \; E)$ of $(n-m)$-dimensional subspaces of $E$, and $G\cdot
\Omega^+(c) \subset \mathrm{Hom}(E, \: F)$ is the open subset of homomorphisms of full rank $m$. The projection
\[
G\times_{P_c} \Omega^+(c) \to G/P_c
\]
is the map which associates to a homomorphism $\varphi\in\mathrm{Hom}(E, \: F)$ its kernel $\mathrm{ker}(\varphi ) \in \mathrm{Grass}(n-m, \: E)$.
In the case $n=m+1$, this is a projective space, and the associated method for proving rationality is called the \emph{double bundle method},
which appeared first in \cite{Bo-Ka}. We discuss this in more detail in Chapter 2.
\end{example}

The discussion in Example \ref{eDoubleBundle} shows that it will be very convenient to have results for the stable rationality of Grassmannians
$\mathrm{Grass}( k, \: E)/G$ (where $E$ is a representation of the reductive group $G$) analogous to Corollary
\ref{cstabratspecialgroups}. One has

\begin{proposition}\xlabel{pStabRatGrassmannians}
Let $E$ be a representation of $G = \mathrm{SL}_p (\CC )$, $p$ prime. Let $X:=\mathrm{Grass} (k, \: E)$ be the Grassmannian of $k$-dimensional
subspaces of
$E$. Assume:
\begin{itemize}
\item
The kernel $Z$ of the action of $G$ on $\PP (E)$ coincides with the center $\ZZ /p\ZZ$ of $\mathrm{SL}_p (\CC )$ and the action of $G/Z$ on $\PP
(E)$ is almost free. Furthermore, the action of
$G$ on
$E$ is almost free and each element of $Z$ not equal to the identity acts homothetically as multiplication by a primitive $p$th root of unity.
\item
$k \le \dim E - \dim G-1$.
\item
$p$ does not divide $k$.
\end{itemize}
Then $X/G$ is stably rational, in fact, $X/G \times \PP^{\dim G +1}$ is rational.
\end{proposition}

\begin{proof}
Let $C_X \subset \Lambda^k (E)$ be the affine cone over $X$ consisting of pure (complety decomposable) $k$-vectors. We will show that under the
assumptions of the proposition, the action of $G$ on $C_X$ is almost free. This will accomplish the proof since $C_X/G$ is generically a torus
bundle over $X/G$ hence Zariski-locally trivial; and the group $G=\mathrm{SL}_p (\CC )$ is special.\\
Let $e_1\wedge e_2 \wedge \dots \wedge e_k$ be a general $k$-vector in $\Lambda^k (E)$. Since $k \le \dim E - \dim G-1$ and, in $\PP (E)$, $\dim
(G\cdot [e_1]) = \dim G$ since $Z$ is finite and $G/Z$ acts almost freely on $\PP (E)$, the $k-1$-dimensional projective linear subspace spanned
by $e_1, \dots , e_k$ in $\PP (E)$ will intersect the $\dim E -1 - \dim G$ codimensional orbit $G\cdot [e_1]$ only in $[e_1]$. Hence, if an element
$g\in G$ stabilizes $e_1\wedge \dots \wedge e_k$, it must lie in $Z$. Thus $g \cdot (e_1 \wedge \dots \wedge e_k) = \zeta^{k} (e_1 \wedge \dots
\wedge e_k)$ for a primitive $p$-th root of unity $\zeta$ if $g\neq 1$. But since $p$ does not divide $k$, the case $g\neq 1$ cannot occur.
\end{proof}

As an application we prove the following result which had not been obtained by different techniques so far.

\begin{theorem}\xlabel{tRationalityV34}
The moduli space $\PP (\mathrm{Sym}^{34} (\CC^3)^{\vee }) /\mathrm{SL}_3 (\CC )$ of plane curves of degree $34$ is rational.
\end{theorem}

\begin{proof}
As usual, $V(a,\: b)$ denotes the irreducible $\mathrm{SL}_3 (\CC )$-module whose highest weight has numerical labels $a,\: b$ (we choose the
standard diagonal torus and Borel subgroup of upper triangular matrices for definiteness). Then
\[
V(0, \: 34) \subset \mathrm{Hom} (V(14,\: 1), \: V(0,\: 21))\, ,
\]
and $\dim V(14,\: 1) = 255$, $\dim V(0,\: 21)=253$, so we get a map
\[
\pi\, :\, \PP (V(0,\: 34)) \dasharrow \mathrm{Grass}(2,\: V(14,\: 1))\, 
\]
$\dim \PP (V(0,\: 34))= 629$ and $\dim \mathrm{Grass}(2,\: V(14,\: 1)) = 506$. Moreover, Proposition \ref{pStabRatGrassmannians} and its proof 
show that $\mathrm{Grass}(2,\: V(14,\: 1))/\mathrm{SL}_3 (\CC) \times \PP^9$ is rational, and the action of $\mathrm{PGL}_3 (\CC ) = \mathrm{SL}_3
(\CC ) / Z$, where $Z$ is the center of $\mathrm{SL}_3 (\CC )$, is almost free on $\mathrm{Grass}(2,\: V(14,\: 1))$. Moreover, let $\mathcal{O}_P (1)$ 
be the $\mathrm{SL}_3 (\CC)$-linearized line bundle induced by the Pl\"ucker embedding on $\mathrm{Grass}(2,\: V(14,\: 1))$:
\[
\mathrm{Grass}(2,\: V(14,\: 1)) \subset \PP (\Lambda^2 (V(14,\: 1)))\, .
\]
If we choose on $\PP (V(0,\: 34)) \times \mathrm{Grass}(2,\: V(14,\: 1))$ the bundle $\mathcal{L}:=\mathcal{O}(1)\boxtimes
\mathcal{O}_P (2)$, all the assumptions of Theorem \ref{tMainTrick} except the dominance of $\pi$ have been checked (compare also Example
\ref{eDoubleBundle}). The latter dominance follows from an explicit computer calculation, where one  has to check that for a random element $x_0$
in $V(0,\: 34)$ the corresponding homomorphism in $\mathrm{Hom} (V(14,\: 1), \: V(0, \: 21))$ has full rank, and the fibre of $\pi$ over $\pi
([x_0])$ has the expected dimension $\dim\PP (V(0,\: 34)) - \dim \mathrm{Grass} (2,\: V(14,\: 1) )$.
\end{proof}

\begin{remark}\xlabel{rNoAlternatives}
As far as we can see, the rationality of $\PP (V(0,\: 34))/\mathrm{SL}_3 (\CC )$ cannot be obtained by direct application of the double bundle
method, i.e. by applying Theorem \ref{tMainTrick} in the case discussed in Example \ref{eDoubleBundle} with base of the projection a projective
space. In fact, a computer search yields that the inclusion $V(0, \: 34) \subset \mathrm{Hom}(V(30, 0), \: V(0,4) \oplus V(5,\: 9))$ is the only
candidate to be taken into consideration for dimension reasons: $\dim V(30, 0) = \dim ( V(0,\: 4) \oplus V(5, \: 9) ) +1$ and $\dim \PP (V (0,\:
34)) > \dim \PP (V(30, \: 0))$. However, on $\PP (V(0, \: 34)) \times \PP (V (30, \: 0))$ there does not exist a $\mathrm{PGL}_3 (\CC
)$-linearized line bundle cutting out $\mathcal{O}(1)$ on the fibres of the projection to $\PP (V(30, 0))$; for such a line bundle would have to
be of the form $\mathcal{O}(1) \boxtimes \mathcal{O}(k)$, $k\in\ZZ $, and none of these is $\mathrm{PGL}_3 (\CC )$-linearized: since
$\mathcal{O}\boxtimes \mathcal{O}(1)$ is $\mathrm{PGL}_3 (\CC )$-linearized it would follow that the $\mathrm{SL}_3 (\CC )$ action on $H^0(\PP
(V(0, \: 34)), \: \mathcal{O}(1))\simeq V(34, \: 0)$ factors through $\mathrm{PGL}_3 (\CC )$ which is not the case.
\end{remark}

\section[Overview of some further topics]{Overview of some further topics}

Here we give a brief description of additional topics connected with the rationality problem which are too important to be omitted altogether, but are
outside the focus of the present text.\\
The first concerns \emph{cohomological obstructions to rationality}. One of the first examples was given by Artin and Mumford \cite{A-M} who showed

\begin{proposition}\xlabel{pArtinMumford}
The torsion subgroup $T\subset H^3 (X, \: \ZZ )$ is a birational invariant of a smooth projective variety $X$. In particular, $T=0$ if $X$ is rational.
\end{proposition}

They use this criterion to construct unirational irrational threefolds $X$, in fact $X$ with $2$-torsion in $H^3(X, \: \ZZ )$. Later, David Saltman
\cite{Sa} proved that there are invariant function fields $\CC (X )$ (where $X=V/G$) of the action of a finite group $G$ which are not purely transcendental
over the ground field $\CC$ using as invariant the \emph{unramified Brauer group} $\mathrm{Br}_{\mathrm{nr}} (\CC (X) /\CC )$ which can be shown to equal
the cohomological Brauer group $\mathrm{Br} (\tilde{X}) = H^2_{\mathrm{\acute{e}t}} (\tilde{X}, \mathbb{G}_m )$ of a smooth projective model $\tilde{X}$ of
$\CC (X)$. Here $\mathbb{G}_m$ denotes the sheaf (for the \'{e}tale site on any scheme) defined by the standard multiplicative group scheme. Moreover, in
the sequel ${\mu}_n$ will, as usual, denote the subsheaf of $\mathbb{G}_m$ defined by $\mu_n (U) =$ group of $n$th roots of $1$ in $\Gamma (U, \:
\mathcal{O}_U)$. The unramified point of view was developed further in \cite{Bogo3}.\\
The notion of unramified cohomology generalizes the previous two examples. A particular feature of the unramified view point is that it bypasses the need to
construct a smooth projective model for a given variety $X$, working directly with the function field of $X$, or rather with all smooth projective models
of $X$ at once. Below, $k$ is an algebraically closed field of any characteristic.

\begin{definition}\xlabel{dUnramifiedCohomology}
Let $X$ be a variety over $k$ and $n>0$ an integer prime to $\mathrm{char}(k)$. The $i$-th unramified cohomology group of $X$ with coefficients in
$\mu_n^{\otimes j}$ is by definition
\begin{gather*}
H^i_{\mathrm{nr}} (k (X)/k, \: \mu_n^{\otimes j}):= 
\bigcap_{A\in \mathrm{DVR}(k(X))} \left( \mathrm{im} \left( H^i_{\mathrm{\acute{e}t}} (A, \: \mu_n^{\otimes j}) \to  H^i_{\mathrm{\acute{e}t}} (k (X) , \:
\mu_n^{\otimes j}) \right) \right)
\end{gather*}
where $A$ runs over all rank one discrete valuation rings $k \subset A \subset k(X)$ with quotient field $k(X)$. The cohomology groups are to be interpreted
as \'{e}tale cohomology
\begin{gather*}
H^i_{\mathrm{\acute{e}t}} (A, \: \mu_n^{\otimes j}) := H^i_{\mathrm{\acute{e}t}} (\mathrm{Spec} (A) , \: \mu_n^{\otimes j}), \\
H^i_{\mathrm{\acute{e}t}} (k (X) , \: \mu_n^{\otimes j}) := H^i_{\mathrm{\acute{e}t}} (\mathrm{Spec}(k(X)) , \: \mu_n^{\otimes j})\, .
\end{gather*}
\end{definition} 

It is known that if $k(X)$ and $k(Y)$ are stably isomorphic over $k$, then 
\[
H^i_{\mathrm{nr}} (k (X)/k, \: \mu_n^{\otimes j})\simeq H^i_{\mathrm{nr}} (k (Y)/k, \: \mu_n^{\otimes j})
\]
and in particular, the higher unramified cohomology groups are trivial if $X$ is stably rational (see \cite{CT95}).

\

Clearly, if $G=\mathrm{Gal}(k(X)^s/k(X))$ is the absolute Galois group of $K:=k(X)$ 
\[
H^i_{\mathrm{\acute{e}t}} (\mathrm{Spec}(k(X)) , \: \mu_n^{\otimes j}) = H^i (G, \: \mu_n^{\otimes j})= H^i (K, \: \mu_n^{\otimes j})\, ,
\]
the latter being a Galois cohomology group \cite{Se97} ($\mu_n$ the group of $n$th roots of $1$ in $K^s$). There is an alternative description of unramified
cohomology in terms of residue maps in Galois cohomology which is often useful. We would like to be as concrete as possible, so recall first that given a
profinite group $G$ and a discrete $G$-module $M$ on which $G$ acts continuously, and denoting by $C^n (G, \: M)$ the set of all continuous maps from $G^n$
to
$M$, we define the cohomology groups $H^q (G, \: M)$ as the cohomology of the complex $C^{\cdot } (G, \: M)$ with differential
\begin{align*}
d\, :\, C^n(G, \: M) \to C^{n+1} (G, \: M) , & \\
(df) (g_1, \dots , g_{n+1}) &= g_1\cdot f(g_2, \dots , g_{n+1})\\
                            & +\sum_{i=1}^n (-1)^i f(g_1, \dots ,g_ig_{i+1}, \dots, g_{n+1})\\ 
                            & + (-1)^{n+1} f(g_1, \dots , g_n)
\end{align*}
and this can be reduced to the finite group case since 
\[
H^q (G, \: M) = \lim_{\rightarrow} H^q (G/U, \: M^U)\, ,
\]
the limit taken over all open normal subgroups $U$ in $G$. We recall from \cite{Se97}, Appendix to Chapter II, the following:
\begin{proposition}\xlabel{pSerre}
If $G$ is a profinite group, $N$ a closed normal subgroup of $G$, $\Gamma$ the quotient $G/N$, and $M$ a discrete $G$-module with trivial action of $N$,
then one has exact sequences for all $i\ge 0$
\[
\begin{CD}
0 \to H^i (\Gamma , \: M) @>\pi >> H^i (G, \: M) @>r>> H^{i-1}(\Gamma, \: \mathrm{Hom}(N, \: M)) @>>> 0
\end{CD}
\]
provided the following two assumptions hold:
\begin{itemize}
\item[(a)]
The extension
\[
\begin{CD}
1 @>>> N @>>> G @>>> \Gamma @>>> 1
\end{CD}
\]
splits.
\item[(b)]
$H^i (N, \: M)=0$ for all $i>1$.
\end{itemize}
\end{proposition}

Here $\pi$ is induced through the map $G \to \Gamma$ by functoriality and $r$ is the residue map which has an explicit description as follows: an element
$\alpha\in H^i (G, \: M)$ can be represented by a cocycle $f=f(g_1, \dots , g_i)\in C^i (G, \: M)$ which is normalized (i.e. equal to $0$ if one $g_j$ is
$1$) and which only depends  on $g_1$ and the classes $\gamma_2, \dots , \gamma_i$ of $g_2, \dots , g_i$ in $\Gamma$. If then $\gamma_2, \dots , \gamma_{i}$
are elements in $\Gamma$, one defines $r(f)(\gamma_2, \dots , \gamma_{i})$ to be the element of $\mathrm{Hom}(N, \: M)$ (continuous homomorphisms of $N$ to
$M$) given by
\[
r(f)(\gamma_2, \dots , \gamma_{i})(n) := f (n, \: g_2, \dots , g_i), \: n\in N\, .
\]
The $(i-1)$-cochain $r(f)$ is checked to be an $(i-1)$-cocycle of $\Gamma$ with values in $\mathrm{Hom}(N, \: M)$, and its class $r(\alpha )$ in
$H^{i-1}(\Gamma, \: \mathrm{Hom}(N, \: M))$ is independent of $f$. The proof uses the spectral sequence of group extensions; item (b) is used to reduce the
spectral sequence to a long exact sequence, and (a) to split the long exact sequence into short exact ones of the type given. Details may be found in
\cite{Se97}.\\
Reverting to our original set-up we have the field extension $K/k=k(X)/k$ and a discrete valuation ring $A\subset K$ with $k\subset A$ and such
that
$K$ is the field of fractions of $A$. To $A$ there is an associated completion $\hat{K}_A$. Assume $\hat{K}_A$ of residue characteristic $0$ for
simplicity. The absolute Galois group of $\hat{K}_A$ splits as $\hat{\mathbb{Z}}\oplus G_A$ with $G_A$ the absolute Galois group of the residue field of
$A$. Since $\hat{\mathbb{Z}}$ has cohomological dimension $1$ we can apply Proposition \ref{pSerre} to obtain a map
\[
\varrho_A \, :\, H^i ( K, \: \mu_n^{\otimes j}) \to H^i (\hat{\mathbb{Z}}\oplus G_A, \: \mu_n^{\otimes j}) \to H^{i-1}(G_A, \: \mu_n^{\otimes {j-1}})\, ,
\]
where the first map is restriction and the second one the residue map of Proposition \ref{pSerre}. Then one has
\begin{gather*}
H^i_{\mathrm{nr}}(k(X)/k, \: \mu_n^{\otimes j}) = \bigcap_{A \in \mathrm{DVR}(k(X))} \left( \mathrm{ker} (\varrho_A ) \right)
\end{gather*}

(cf. \cite{CT95} or \cite{CT-O}, \S 1 for the proof). This is a purely Galois cohomological description of unramified cohomology.

\

To connect the notion of unramified cohomology with the classical work of Artin-Mumford and Saltman, we list a few results when unramified cohomology has
been computed.

\begin{itemize}
\item
Let $X$ be a smooth projective variety over a field $k$, algebraically closed of characteristic $0$. Then
\[
H^1_{\mathrm{nr}} (k(X)/k, \mu_n)\simeq (\ZZ /n\ZZ )^{\oplus 2q} \oplus _n \mathrm{NS}(X)
\]
$q=\dim H^1 (X, \:\mathcal{O}_X)$ is the dimension of the Picard variety of $X$ and $_n\mathrm{NS}(X)$ is $n$-torsion in the N\'{e}ron-Severi
group of $X$. Cf. \cite{CT95}, Prop. 4.2.1.
\item
Let $X$ be as in the previous example, suppose furthermore that $X$ is unirational and let $n$ be a power of a prime number $l$. Then
\[
H^2_{\mathrm{nr}} (k(X)/k, \mu_n)\simeq _n H^3(X, \ZZ_l)
\]
(where $H^3(X, \: \ZZ_l)$ is the third \'{e}tale cohomology group of $X$ with $\ZZ_l$-coefficients). Furthermore, $H^2_{\mathrm{nr}} (k(X)/k, \mu_n)\simeq
_n \mathrm{Br}(X)$ which explains the relation to the examples of Artin-Mumford and Saltman. Cf. \cite{CT95}, Prop. 4.2.3.
\end{itemize}

Much less is known about higher unramified cohomology groups\\ $H^i_{\mathrm{nr}}(k(X)/k, \: \mu_n^{\otimes j})$, $i\ge 3$, but see \cite{CT-O}, \cite{Sa3}
and
\cite{Pey}, \cite{Mer} for some computations and uses of unramified $H^3$. Let us mention at this point a problem which seems particularly attractive and
which is apparently unsolved.

\begin{problem}\xlabel{pUnramifiedPGL8}
Since for any divisor $n$ of $420$ and any almost free representation $V$ of $\mathrm{PGL}_n (\CC )$, the field $\CC (V)^{\mathrm{PGL}_n (\CC )}$ is stably
rational (see \cite{CT-S}, Prop. 4.17 and references there), it would be very interesting to compute some higher unramified cohomology groups of $\CC
(V)^{\mathrm{PGL}_8 (\CC )}$, for some almost free representation $V$ of $\mathrm{PGL}_8 (\CC )$, or at least to detect nontrivial elements in some such
group (higher should probably mean here degree at least $4$). This would give examples of $\mathrm{PGL}_8 (\CC )$-quotients which are not stably rational.
\end{problem}

We conclude by some remarks how the rationality problem changes character if we allow our ground field $k$ to be non-closed or of positive characteristic.
First, Merkurjev \cite{Mer} has shown that over nonclosed fields $k$ there exist examples of connected simply connected semi-simple groups $G$ with almost
free action on a linear representation $V$ such that $k(V)^G$ is non-rational (even not stably rational). On the other hand, over $k=\bar{\mathbb{F}}_p$, 
the recent article \cite{BPT-2} proves stable rationality for many quotients $V/G$ where $G$ is a finite group $G$ of Lie type, and $V$ a faithful
representation of $G$ over $k=\bar{\mathbb{F}}_p$.


\chapter[Techniques and recent results]{Techniques for proving rationality and some recent results for moduli spaces of plane curves}

\section{Introduction}

In this chapter we collect, in all brevity, the main methods for approaching the rationality problem (Problem \ref{problemrat}). In each case we list
some illustrative results obtained by the respective method as a guide to the literature, or develop applications in the text itself.\\
In the last section we summarize some more recent results on the rationality problem for moduli spaces of plane curves of fixed degree (under
projectivities).

\section{Methods}

\subsection{Slice method and 2-form trick}

Given a linear algebraic group $G$ and $G$-variety $X$, the study of birational properties of the quotient $X/G$ can be reduced to a smaller variety
and smaller group $H < G$ in the following way.

\begin{definition}\xlabel{dGHsection}
A $(G, \: H)$-section of $X$ is an irreducible $H$-stable subvariety $Y$ in X whose translates by elements of $G$ are dense in $X$, and with the
property that if $g\in G$ carries two points in $Y$ into one another, then $g$ is already in $H$.
\end{definition}

Then $Y/H$ is birational to $X/G$ as can be seen by restricting rational functions on $X$ to $Y$ which induces an isomorphism of invariant function
fields. The \emph{slice method}, consisting in finding a $(G, \: H)$-section as above, is not so much a direct method for proving rationality as
rather a preliminary or intermediate simplification step applied in the course of the study of the birational properties of a given space. As such it
corresponds to the simple idea of \emph{reduction to normal form}.\\
The so called \emph{2-form trick} (see \cite{Shep}, Prop. 8) has already been
mentioned above, and is contained in Theorem
\ref{tMainTrick} in conjunction with Example \ref{eTwoForms}. We phrase it here again in more explicit form for reference.

\begin{theorem}\xlabel{tTwoFormTrick}
Let $E$ be a finite dimensional representation of odd dimension of a reductive group $G$, and let $V$ be a subrepresentation of $\Lambda^2
(E)^{\vee}$. Let $Z$ be the kernel of the action of $G$ on $\PP (E)$, and suppose $Z$ acts trivially on $\PP (V)$. Assume that the action of $G/Z$ on $\PP
(E)$ is almost free and that there exists a $G/Z$-linearized line bundle 
$\mathcal{L}$ on
$\PP (V) \times \PP (E)$ such that $\mathcal{L}$ cuts
out
$\mathcal{O}(1)$ on the fibres of the projection to $\PP (E)$.\\ Suppose that for some $v_0\in V$ the associated $2$-form in $\Lambda^2 (E)^{\vee }$, viewed
as a skew-symmetric map $E \to E^{\vee}$, has maximal rank
$\dim E -1$. Then the rational map
\[
\varphi \, :\, \PP (V) \dasharrow \PP (E) \, ,
\]
associating to a $2$-form its kernel, is well-defined, and if $\dim V > \dim E$, $\dim \varphi^{-1}(\varphi ([v_0])) = \dim V -\dim E$, then $\varphi$
is dominant. Hence, if $\PP (E)/G$ is stably rational of level $\le \dim V -\dim E$, then $\PP (V)/G$ is rational.
\end{theorem}

This method was used in \cite{Shep} to prove the rationality of the moduli spaces $\PP (\mathrm{Sym}^d (\CC^3)^{\vee})/\mathrm{SL}_3 (\CC )$ of plane
curves of degrees $d\equiv 1$ (mod $4$).\\
To give an illustration of both the slice method and $2$-form trick we will study the invariant function field
\[
\mathbb{K}_d =\CC \left( \mathrm{Sym}^2 (\CC^d)^{\vee } \otimes \CC^3 \right)^{\mathrm{SL}_d (\CC) \times \mathrm{SL}_3 (\CC ) \times \CC^{\ast}} \, .
\]
Recall that a \emph{theta-characteristic} $\theta$ on a smooth plane curve $C$ of degree $d$ is a line bundle which is a square root of the canonical
line bundle $\omega_C$. Then the above function field is the field of rational functions on the moduli space of pairs $(C, \: \theta )$ where $C$ is a
smooth plane curve of degree $d$ as above and $\theta$ a theta-characteristic on $C$ with $h^0 (C, \: \theta )=0$, see \cite{Beau00}. So $C$ is the
discriminant curve of a net of quadrics in $\PP^{d-1}$; the above invariant function field may also be interpreted as the moduli space of nets of
quadrics in $\PP^{d-1}$.\\
At the 2008 Geometric Invariant Theory Conference in G\"ottingen, F. Catanese asked whether $\mathbb{K}_d$ was rational; this may in general be quite
tricky to decide, in particular for $d$ even (note that $\mathbb{K}_4$ is the function field of $\mathfrak{M}_3$, and the proof of rationality of
$\mathfrak{M}_3$ presented great difficulties, see \cite{Kat92/2}, \cite{Kat96}). We assume $d$ odd in the sequel and show how the problem may be reduced to
the question of rationality for a simpler invariant function field (this approach suggests that the problems of rationality of $\mathbb{K}_d$ for $d$ odd
and even are interrelated, and there could be an inductive procedure for proving rationality for all $d$). For $d=5$ this approach was worked out in
\cite{Kat92/1}.\\ Note that canonically $\CC^3 \simeq \Lambda^2 (\CC^3)^{\vee }$, so that
\begin{gather*}
\mathrm{Sym}^2 (\CC^d)^{\vee } \otimes \CC^3 \simeq \mathrm{Sym}^2 (\CC^d)^{\vee } \otimes \Lambda^2 (\CC^3)^{\vee } \subset \Lambda^2
((\CC^d)^{\vee}\otimes (\CC^3)^{\vee})
\end{gather*}
so that by the 2-form trick we may obtain a rational map
\[
\varphi\, :\, \mathrm{Sym}^2 (\CC^d)^{\vee } \otimes \Lambda^2(\CC^3)^{\vee } \dasharrow \PP (\CC^d \otimes \CC^3)\, .
\]
We need
\begin{lemma}\xlabel{lTechnicalTheta}
The inclusion 
\[
\iota\, :\, \mathrm{Sym}^2 (\CC^d)^{\vee } \otimes \Lambda^2 (\CC^3)^{\vee } \subset \Lambda^2
((\CC^d)^{\vee}\otimes (\CC^3)^{\vee})
\]
is given in terms of coordinates $x_1, \dots , x_d$ in $(\CC^d)^{\vee }$ and $y_1,\: y_2, \: y_3$ in $(\CC^3)^{\vee }$ by
\[
\iota ( (x_ix_j) \otimes (y_k \wedge y_l)) = (x_i \otimes y_k) \wedge (x_j \otimes y_l) + (x_j \otimes y_k)\wedge (x_i \otimes y_l)\, .
\]
The map $\varphi\, :\, \mathrm{Sym}^2 (\CC^d)^{\vee } \otimes \Lambda^2(\CC^3)^{\vee } \dasharrow \PP (\CC^d \otimes \CC^3)$ is dominant for odd $d\ge 5$. 
\end{lemma}
\begin{proof}
We just have to check the well-definedness and dominance of $\varphi$. If $d=5$ then the element
\begin{gather*}
\omega = (x_5^2-2x_1x_2)\otimes (y_2\wedge y_3) + (x_1^2+x_3^2+x_4^2)\otimes (y_1\wedge y_3)\\ + (2x_4x_5-2x_2x_3) \otimes (y_1 \wedge y_2)
\end{gather*}
when viewed as an element of $\Lambda^2 ((\CC^5)^{\vee }\otimes (\CC^3)^{\vee })$ or alternatively an antisymmetric $15\times 15$-matrix, has rank
exactly $14$. Thus, since $d-5$ is even, it suffices to indicate an element in
\[
 \mathrm{Sym}^2 (\CC^2)^{\vee } \otimes \Lambda^2 (\CC^3)^{\vee } \subset \Lambda^2 ((\CC^2)^{\vee }\otimes (\CC^3)^{\vee })
\]
whose associated $6\times 6$ antisymmetric matrix has maximal rank: we decompose $\CC^{d-5}\otimes \CC^3 = (\CC^2\otimes\CC^3) \oplus \dots \oplus
(\CC^2\otimes\CC^3)$ ($(d-5)/2$ times) and consider an antisymmetric matrix of block diagonal form with one $15\times 15$ block and $(d-5)/2$ blocks of
size $6\times 6$. An element of the required form is
\[
\pi_j= x_j^2\otimes (y_1\wedge y_2)  + (x_jx_{j+1})\otimes (y_1\wedge y_3) + x_{j+1}^2 \otimes (y_2\wedge y_3)
\]
(here $j$ runs over the even numbers between $6$ and $d-1$). Thus we have found an element in the image of $\iota$ with one-dimensional kernel for
every $d$, namely
\[
\kappa = \omega + \sum_j \pi_j \, .
\] 
Thus $\varphi$ is well-defined for $d\ge 5$, $d$ odd. Moreover, the kernel of $\iota (\kappa )$ is spanned exactly by the matrix in $\CC^d \otimes
\CC^3$ which in terms of coordinates $e_1,\dots , e_d$ in $\CC^d$ and $f_1, \: f_2, \: f_3$ in $\CC^3$ dual to the $x_i$ and $y_j$ is 
\[
m= e_1\otimes f_1 + e_2\otimes f_2 + e_3\otimes f_3 
\]
(note that it suffices to check this for $d=5$, since the $\pi_j$ vanish on $m$ by construction). Since $\mathrm{SL}_d (\CC ) \times \mathrm{SL}_3 (\CC
)$ has a dense orbit on $\PP (\CC^d \otimes \CC^3)$ (the matrices of maximal rank), and $m$ has maximal rank, this checks dominance of $\varphi$.
\end{proof}

Thus we see that if we put
\[
L:= \overline{ \varphi^{-1} ([m] )} \subset \mathrm{Sym}^2 (\CC^d)^{\vee }\otimes \CC^3
\]
then $L$ is a linear subspace which is a $(\mathrm{SL}_d (\CC ) \times \mathrm{SL}_3 (\CC ), \: H)$-section of $\mathrm{Sym}^2 (\CC^d)^{\vee }\otimes
\CC^3$ where $H$ is the stabilizer of $[m] \in \PP (\CC^d \otimes \CC^3)$ in $\mathrm{SL}_d (\CC ) \times \mathrm{SL}_3 (\CC )$. Moreover,
\[
\mathbb{K}_d \simeq \CC (L)^{H\times\CC^{\ast}}\, .
\]
We would like to describe the $H$-representation $L$ more explicitly. We note that $\mathrm{SL}_d (\CC )\times\mathrm{SL}_3 (\CC )$ acts on $\CC^d
\otimes \CC^3$, viewed as $d\times 3$-matrices, as
\[
(A, \: B)\cdot M = AMB^t, \quad (A, \: B) \in \mathrm{SL}_d (\CC )\times\mathrm{SL}_3 (\CC ), \: M \in \mathrm{Mat}(d\times 3, \: \CC )
\]
and consequently
\begin{gather*}
H = \left\{ \left( \left( \left(
\begin{array}{cc} \lambda^{d-3} s & 0 \\ \ast & \lambda^{-3} S  \end{array}
\right)^t \right)^{-1} , \: s \right) \, \mid \, S \in \mathrm{SL}_{d-3} (\CC ), \: s\in \mathrm{SL}_3 (\CC ), \: \lambda\in\CC^{\ast }   \right\}\, .
\end{gather*}

We introduce some further notation. We denote by
\[
P = \left\{ \left( \begin{array}{cc} s & 0 \\ \ast & S \end{array} \right) \, \mid \, S \in \mathrm{GL}_{d-3} (\CC ), \: s\in \mathrm{GL}_3 (\CC )
\right\} \subset \mathrm{GL}_d (\CC )
\] 
the indicated parabolic subgroup of $\mathrm{GL}_d (\CC )$ and put
\[
P' = \left\{ \left( \begin{array}{cc} s & 0 \\ \ast & S \end{array} \right) \, \mid \, S \in \mathrm{SL}_{d-3} (\CC ), \: s\in \mathrm{SL}_3 (\CC )
\right\} \, .
\]
$P'$ is a subgroup of $H$ in the natural way; we will investigate the structure of $L$ as $P'$-module first and afterwards do the bookkeeping for the
various torus actions. Associated to the standard representation of $\mathrm{GL}_d (\CC )$ on $\CC^d$ and the parabolic $P$ we have the $P$-invariant
subspace $F$ below and complement $E$:
\[
F:= \langle e_4, \dots , e_d \rangle \subset \CC^d, \: E:= \langle e_1, \: e_2, \: e_3 \rangle \subset \CC^d \, ,
\]
so that we have a filtration
\begin{gather*}
\mathrm{Sym}^3 (F) \subset \mathrm{Sym}^3 (F) \oplus \mathrm{Sym}^2 (F)\otimes E \subset \\ \mathrm{Sym}^3 (F) \oplus \mathrm{Sym}^2 (F)\otimes E \oplus
F \otimes \mathrm{Sym}^2 (E) \subset \mathrm{Sym}^3 (\CC^d)\, .
\end{gather*}
Then we claim
\begin{lemma}\xlabel{lThetaMain}
There is an isomorphism of $P'$-modules (for $d\ge 5$ odd)
\[
L \simeq  \mathrm{Sym}^3 (\CC^d) /\mathrm{Sym}^3 (F) \, .
\]
\end{lemma}

\begin{proof}
We first remark that the dimensions are right: in fact
\begin{align*}
\dim L &= \dim \mathrm{Sym}^2 (\CC^d)^{\vee}\otimes \CC^3 - \dim \PP (\CC^d \otimes \CC^3)\\
       &= 3 {d+1 \choose 2} - (3d-1) = \frac{3}{2} d^2 -\frac{3}{2} d +1
\end{align*}
whereas 
\begin{align*}
\dim (\mathrm{Sym}^3 (\CC^d) /\mathrm{Sym}^3 (F) ) &= \dim (\mathrm{Sym}^3 (\CC^d)) - \dim (\mathrm{Sym}^3 (F))\\
    &= { d+2 \choose 3} - { d-1 \choose 3}\\ &= \frac{(d+2)(d+1)d}{6} - \frac{(d-1)(d-2)(d-3)}{6}\\
    &= \frac{3}{2} d^2 -\frac{3}{2} d +1 \, .
\end{align*}
We will construct a $P'$-isomorphism $\mathrm{Sym}^3 (\CC^d) /\mathrm{Sym}^3 (F) \to L$. The representation of the stabilizer
group $H\subset \mathrm{SL}_d (\CC ) \times \mathrm{SL}_3 (\CC )$ in $\mathrm{Sym}^2 (\CC^d)^{\vee }\otimes \Lambda^2 (\CC^3)^{\vee }$ induces a
representation of the subgroup $P' \subset H$ which is isomorphic to 
\begin{gather*}
P' \to \mathrm{Aut}\left( \mathrm{Sym}^2 (\CC^d)\otimes \CC^3 \right)
\end{gather*} 
where the action of $P'$ on $\mathrm{Sym}^2 (\CC^d)$ is obtained by restricting the usual action of $\mathrm{SL}_d (\CC )$ to the subgroup $P'$, and
the action of $P'$ on $\CC^3$ is given by
\[
\left( \begin{array}{cc} s & 0 \\ \ast & S \end{array} \right) \cdot v = s\cdot v\, .
\]
We identify the representation $L$ in this picture by showing that there is a unique subspace of dimension
\[
\frac{3}{2} d^2 -\frac{3}{2} d +1
\]
in $\mathrm{Sym}^2 (\CC^d)\otimes \CC^3$ which is invariant under the semisimple subgroup $\mathrm{SL}_{d-3}(\CC ) \times \mathrm{SL}_3 (\CC ) \subset
P'$. One has
\begin{align*}
\mathrm{Sym}^2 (\CC^d) \otimes \CC^3 &\simeq \mathrm{Sym}^2 (F \oplus E) \otimes E\\
 &\simeq \left( (\mathrm{Sym}^2 (F)) \oplus (F\otimes E) \oplus (\mathrm{Sym}^2 (E)) \right) \otimes E\\
 &\simeq (\mathrm{Sym}^2 (F)\otimes E) \oplus (F \otimes \Lambda^2 (E)) \oplus (F \otimes\mathrm{Sym}^2(E))\\
     &\oplus \mathrm{Sym}^3 (E) \oplus \Sigma^{2,\, 1}( E)
\end{align*}
where the last two lines give the decomposition into irreducible $\mathrm{SL}_{d-3}(\CC ) \times \mathrm{SL}_3 (\CC )$-modules. The dimensions of these,
listed in the order in which they occur in the last two lines of the previous formula, are
\[
3\cdot {d-2 \choose 2}, \; 3(d-3), \; 6 (d-3), \; 10, \; 8\, .
\]
Now we always have
\[
\frac{3}{2} d^2 -\frac{3}{2} d +1-(3(d-3)+ 6 (d-3) + 10 + 8) = \frac{3}{2}d^2 -\frac{21}{2} d +10 > 0
\]
as soon as $d\ge 7$. Thus then $(\mathrm{Sym}^2 (F)\otimes E) \subset L$ and
\[
\dim L - \dim (\mathrm{Sym}^2 (F)\otimes E) = 6d -8\, . 
\]
Certainly, $3(d-3) +10+8 = 3d + 9 < 6d -8$ for $d\ge 7$, so we find
\[
L =(\mathrm{Sym}^2 (F)\otimes E) \oplus (F \otimes\mathrm{Sym}^2(E)) \oplus \mathrm{Sym}^3 (E)\, .
\] 
This is also true for $d=5$: here $\dim L= 31$ and the dimensions of the previous representations are $9, \: 6, \: 12, \: 10, \: 8$. Experimenting a
little shows that we have to take again the $9$, $12$ and $10$ dimensional representations to get $31$. Now
\[
\mathrm{Sym}^2 (\CC^d) \otimes \CC^3 \simeq (\mathrm{Sym}^2 (\CC^d) \otimes \CC^d)/ (\mathrm{Sym}^2 (\CC^d) \otimes F)
\]
as $P'$-representations. The composition of the inclusion with the projection:
\[
\mathrm{Sym}^3(\CC^d) \subset \mathrm{Sym}^2(\CC^d)\otimes \CC^d \to (\mathrm{Sym}^2 (\CC^d) \otimes \CC^d)/ (\mathrm{Sym}^2 (\CC^d) \otimes F)
\]
is $P'$-equivariant, and again viewing this as a map of $\mathrm{SL}_{d-3}(\CC )\times \mathrm{SL}_3 (\CC )$-modules or, equivalently, using the
splitting $\CC^d = F\oplus E$, we find that this map induces the desired $P'$-isomorphism
\[
\mathrm{Sym}^3 (\CC^d) / \mathrm{Sym}^3 (F) \to L\, .
\]
Remark that $\mathrm{Sym}^3 (E) \subset \mathrm{Sym}^3 (\CC^d)$ maps to $L$ nontrivially, whence also the copies of $\mathrm{Sym}^2 (F)\otimes E$ and
$F\otimes\mathrm{Sym}^2 (E)$ contained in $\mathrm{Sym}^3 (\CC^d )$ map to $L$ nontrivially by $P'$-invariance.
\end{proof}

We can now easily obtain

\begin{proposition}\xlabel{pThetaIso}
For the field \[\mathbb{K}_d =\CC \left( \mathrm{Sym}^2 (\CC^d)^{\vee } \otimes \CC^3 \right)^{\mathrm{SL}_d (\CC) \times \mathrm{SL}_3 (\CC ) \times
\CC^{\ast}}\] one has the isomorphism
\[
\mathbb{K}_d \simeq \CC( \mathrm{Sym}^3 (\CC^d)/\mathrm{Sym}^3 (F))^P\, .
\]
\end{proposition}

\begin{proof}
We have seen $\mathbb{K}_d \simeq \CC (L)^{H\times \CC^{\ast}}$, $\CC^{\ast}$ acting by homotheties. The group $P$ is generated by the group $P'$ and
the two dimensional torus
\[
T= \left\{  \left( \begin{array}{cc} t_1 \mathrm{Id}_{3} & 0 \\ 0 & t_2 \mathrm{Id}_{d-3} \end{array} \right) \, \mid \, t_1, \: t_2 \in \CC^{\ast}
\right\}\, .
\]
On the other hand, one may view $P'$ as a subgroup of $H\times \CC^{\ast}$ via the assignment
\begin{gather*}
\left(
\begin{array}{cc}  s & 0 \\ \ast &  S  \end{array}
\right)
\mapsto 
\left( \left( \left( \left(
\begin{array}{cc}  s & 0 \\ \ast &  S  \end{array}
\right)^t \right)^{-1} , \: s \right) , \: 1 \right)
\end{gather*}
and there is also a two dimensional torus $\CC^{\ast}\times \CC^{\ast }$ embedded into $H \times \CC^{\ast}$ via
\begin{gather*}
(\lambda ,\: \mu ) \mapsto 
\left( \left( \left( \left(
\begin{array}{cc}  \lambda^{d-3} \mathrm{Id}_3 & 0 \\ \ast &  \lambda^{-3} \mathrm{Id}_{d-3}  \end{array}
\right)^t \right)^{-1} , \: \mathrm{Id}_3 \right) , \: \mu \right)
\end{gather*}
and $H\times \CC^{\ast }$ is generated by $P'$ and $\CC^{\ast}\times \CC^{\ast }$. We also know that there is the isomorphism
\[
\mathrm{Sym}^3 (\CC^d)/ \mathrm{Sym}^3 (F) \to L
\]
of $P'$-modules. Thus to prove the Proposition, it is sufficient to show that under this isomorphism $T$-orbits transform into $\CC^{\ast}\times
\CC^{\ast}$-orbits. This is straightforward to check: we have seen that as $\mathrm{SL}_{d-3} (\CC )\times \mathrm{SL}_3 (\CC)$-representations
\[
L =(\mathrm{Sym}^2 (F)\otimes E) \oplus (F \otimes\mathrm{Sym}^2(E)) \oplus \mathrm{Sym}^3 (E)
\]
and then 
\[
\left( \begin{array}{cc} t_1 \mathrm{Id}_{3} & 0 \\ 0 & t_2 \mathrm{Id}_{d-3} \end{array} \right) \in T
\]
acts as a homothety on $\mathrm{Sym}^{2-i}(F)\otimes \mathrm{Sym}^{i+1}(E)$, $i=0, \: 1, \: 2$, namely as multiplication by $t_2^{2-i}t_1^{i+1}$; but
$(\lambda,\:
\mu)
\in
\CC^{\ast}\times
\CC^{\ast}$ likewise acts via homotheties on the irreducible summands of $L$ (as $\mathrm{SL}_{d-3} (\CC )\times \mathrm{SL}_3 (\CC)$-representation) viewed
as a subspace in $\mathrm{Sym}^2 (\CC^d)^{\vee }\otimes \Lambda^2 (\CC^3)^{\vee}$, namely by multiplication by $\lambda^{Ai+B}\mu$, some $A$ and $B$ in
$\mathbb{Z}$. This concludes the proof of the Proposition.
\end{proof}

As we pointed out above, Proposition \ref{pThetaIso} does not solve the initial rationality problem for plane curves of odd degree with
theta-characteristic, but we thought it useful to record this important reduction step as an illustration of the methods introduced in this section.\\
The decomposition 
\[
L =(\mathrm{Sym}^2 (F)\otimes E) \oplus (F \otimes\mathrm{Sym}^2(E)) \oplus \mathrm{Sym}^3 (E)\, .
\]
suggests that there may be an inductive procedure to reduce rationality of $\mathbb{K}_d$ to rationality of $\mathbb{K}_{d-3}$.

\subsection{Double bundle method}

The main technical point is the so called "no-name lemma".

\begin{lemma}\xlabel{lNoNameLemma}
Let $G$ be a linear algebraic group with an almost free action on a variety $X$. Let $\pi \, :\, \mathcal{E} \to X$ be a $G$-vector bundle of rank
$r$ on $X$. Then one has the following commutative diagram of $G$-varieties
\xycenter{
\mathcal{E} \ar@{-->}[r]^{f} \ar[rd]_{\pi} &X\times\mathbb{A}^r \ar[d]^{\mathrm{pr}_1}\\
&X 
	}
where $G$ acts trivially on $\mathbb{A}^r$, $\mathrm{pr}_1$ is the projection onto $X$, and the rational map $f$ is \emph{birational}.
\end{lemma}

If $X$ embeds $G$-equivariantly in $\PP (V)$, $V$ a $G$-module, $G$ is reductive and $X$ contains stable points of $\PP (V)$, then this is an
immediate application of descent theory and the fact that a vector bundle in the \'{e}tale topology is a vector bundle in the Zariski topology. The
result appears in \cite{Bo-Ka}. A proof without the previous technical restrictions is given in \cite{Ch-G-R}, \S 4.3.\\
The following result (\cite{Bo-Ka}, \cite{Kat89}) is the form in which Lemma \ref{lNoNameLemma} is most often applied since it allows one to extend its
scope to irreducible representations. 

\begin{theorem}\xlabel{tDoubleBundleOriginal}
Let $G$ be a linear algebraic group, and let $U$, $V$ and $W$, $K$ be (finite-dimensional) $G$-representations. Assume that the stabilizer in general
position of $G$ in $U$, $V$ and $K$ is equal to one and the same subgroup $H$ in $G$ which is also assumed to equal the ineffectiveness kernel in these
representations (so that the action of $G/H$ on $U$, $V$, $K$ is almost free).\\
The relations $\dim U -\dim W =1$ and $\dim V - \dim U > \dim K$ are required to hold.\\
Suppose moreover that there is a $G$-equivariant bilinear map 
\[
\psi \, :\, V \times U \to W
\]
and a point $(x_0, \: y_0) \in V\times U$ with $\psi (x_0, \: y_0)=0$ and $\psi (x_0, \: U) =W$, $\psi (V, \: y_0) =W$.\\
Then if $K/G$ is rational, the same holds for $\PP (V)/G$.
\end{theorem}

\begin{proof}
We abbreviate $\Gamma := G/H$ and let $\mathrm{pr}_U$ and $\mathrm{pr}_V$ be the projections of $V\times U$ to $U$ and $V$. By the genericity
assumption on $\psi$, there is a unique irreducible component $X$ of $\psi^{-1}(0)$ passing through $(x_0, \: y_0)$, and there are non-empty open
$\Gamma$-invariant sets
$V_0\subset V$ resp.
$U_0\subset U$ where $\Gamma$ acts with trivial stabilizer and the fibres $X\cap \mathrm{pr}_V^{-1}(v)$ resp. $X\cap \mathrm{pr}_U^{-1}(u)$ have the
expected dimensions $\dim U -\dim W=1$ resp. $\dim V -\dim W$. Thus
\[
\mathrm{pr}_V^{-1}(V_0) \cap X \to V_0 , \quad \mathrm{pr}_U^{-1}(U_0) \cap X \to U_0
\]
are $\Gamma$-equivariant bundles, and by Lemma \ref{lNoNameLemma} one obtains vector bundles
\[
(\mathrm{pr}_V^{-1}(V_0) \cap X)/\Gamma \to V_0/\Gamma , \quad (\mathrm{pr}_U^{-1}(U_0) \cap X)/\Gamma \to U_0/\Gamma
\]
of rank $1$ and $\dim V -\dim W$ 
and there is still a homothetic $T:=\CC^{\ast}\times \CC^{\ast}$-action on these bundles. By Theorem \ref{Rosenlichtsections}, the action of the torus
$T$ on the respective base spaces of these bundles has a section over which the bundles are trivial; thus we get 
\[
\PP (V)/\Gamma \sim (\PP (U) /\Gamma ) \times \PP^{\dim V-\dim W-1} = (\PP (U) /\Gamma ) \times \PP^{\dim V-\dim U} \, .
\]
On the other hand, one may view $U\oplus K$ as a $\Gamma$-vector bundle over both $U$ and $K$; hence, again by Lemma \ref{lNoNameLemma},
\[
U/\Gamma \times \PP^{\dim K} \sim K/\Gamma \times \PP^{\dim U}\, .
\] 
Since $U/\Gamma$ is certainly stably rationally equivalent to $\PP (U)/\Gamma$ of level at most one, the inequality $\dim V -\dim U > \dim K$ insures
that $\PP (V)/\Gamma$ is rational as $K/\Gamma$ is rational. 
\end{proof}

In \cite{Kat89} this is used to prove the rationality of the moduli spaces $\PP (\mathrm{Sym}^d (\CC^3)^{\vee }) /\mathrm{SL}_3 (\CC )$ of plane curves
of degree $d \equiv 0$ (mod $3$) and $d \ge 210$. A clever inductive procedure is used there to reduce the genericity requirement for the occurring
bilinear maps $\psi$ to a purely numerical condition on the labels of highest weights of irreducible summands in $V$, $U$, $W$. This method is only
applicable if $d$ is large.\\
Likewise, in \cite{Bo-Ka}, the double bundle method is used to prove the rationality of $\PP (\mathrm{Sym}^d \CC^2 )/\mathrm{SL}_2 (\CC )$, the moduli
space of $d$ points in $\PP^1$, if $d$ is even.

\subsection{Method of covariants}

Virtually all the methods for addressing the rationality problem (Problem \ref{problemrat}) are based on introducing some fibration structure over a
stably rational base in the space for which one wants to prove rationality; with the Double Bundle Method, the fibres are linear, but it turns out
that fibrations with nonlinear fibres can also be useful if rationality of the generic fibre of the fibration over the function field of the base can
be proven. The \emph{Method of Covariants} (see \cite{Shep}) accomplishes this by inner linear projection of the generic fibre from a very
singular centre.
 
\begin{definition}\xlabel{dCovariants}
If $V$ and $W$ are $G$-modules for a linear algebraic group $G$, then a covariant $\varphi$ of degree $d$ from $V$ with values in $W$ is simply a
$G$-equivariant polynomial map of degree $d$
\[
\varphi\, :\, V \to W\, .
\]
In other words, $\varphi$ is an element of $\mathrm{Sym}^d (V^{\vee})\otimes W$.
\end{definition}

The method of covariants, phrased in a way that we find quite useful, is contained in the following theorem.

\begin{theorem}\xlabel{tCovariants}
Let $G$ be a connected linear algebraic group the semi-simple part of which is a direct product of groups of type $\mathrm{SL}$ or $\mathrm{Sp}$. Let
$V$ and $W$ be $G$-modules, and suppose that the action of $G$ on $W$ is generically free. Let $Z$ be the ineffectivity kernel of the action of $G$ on
$\PP (W)$, and assume that the action of $\bar{G} := G/Z$ is generically free on $\PP (W)$, and $Z$ acts trivially on $\PP (V)$.\\
Let 
\[
\varphi\, :\, V \to W
\]
be a (non-zero) covariant of degree $d$. Suppose the following assumptions hold:
\begin{itemize}
\item[(a)]
$\PP (W)/G$ is stably rational of level $\le \dim \PP (V) - \dim \PP (W)$.
\item[(b)]
If we view $\varphi$ as a map $\varphi\, :\, \PP (V) \dasharrow \PP (W)$ and denote by $B$ the base \emph{scheme} of $\varphi$, then there is a linear
subspace $L \subset V$ such that $\PP (L)$ is contained in $B$ together with its full infinitesimal neighbourhood of order $(d-2)$, i.e.
\[
\mathcal{I}_{B} \subset \mathcal{I}_{\PP (L)}^{d-1}\, .
\]
Denote by $\pi_L$ the projection $\pi_L \, : \, \PP (V) \dasharrow \PP (V /L)$ away from $\PP (L)$ to $\PP (V/L)$.
\item[(c)]
Consider the diagram
\xycenter{
	\PP (V) \ar@{-->}[r]^{\varphi} \ar@{-->}[d]^{\pi_{L}}& \PP (W)\\
	\PP (V/L)
	}
and assume that one can find a point $[\bar{p}] \in \PP (V/L)$ such that
\[
	\varphi |_{\mathbb{P} (L+\mathbb{C} p )} \colon \mathbb{P} (L +\mathbb{C} p)  \dasharrow  \PP (W)
\]
is dominant.
\end{itemize}
Then $\PP (V) /G$ is rational.
\end{theorem}

\begin{proof}
By Corollary \ref{cspecialandtorus}, the projection $\PP (W) \dasharrow \PP (W)/G$ has a rational section $\sigma$. Remark that property (c) implies
that the generic fibre of $\pi_L$ maps dominantly to $\PP (W)$ under $\varphi$, which means that the generic fibre of $\varphi$ maps dominantly to $\PP 
(V/L)$ under $\pi_L$, too. Note also that the map $\varphi$ becomes linear on a fibre $\PP (L +\CC g)$ because of property (b) and that thus the
generic fibre of $\varphi$ is birationally a vector bundle via $\pi_L$ over the base $\PP (V/L)$.  Thus, if we introduce the graph
\[
	\Gamma = \overline{\{ ([q],[\bar{q}], [f]) \,|\, \pi_L ([q])=[\bar{q}], \varphi ([q])= [f]\} }
	\subset \PP (V) \times \PP\bigl( V /L \bigr) \times \PP (W)
\]
and look at the diagram
\newcommand{\rmpr}{\mathrm{pr}}
\xycenter{
	\Gamma \ar[d]_{\rmpr_{23}} \ar@{<-->}[r]^{1:1}_{\rmpr_{1}}
	& \PP (V) \ar@{-->}[r]
	& \PP (V)/\bar{G} \ar@{-->}[dd]^{\bar{\varphi}}
	\\
	\PP\bigl( V/L \bigr) \times \PP (W) \ar[d]
	\\
	\PP (W) \ar@{-->}[rr]
	& & \PP (W)/\bar{G} \ar@/^20pt/@{-->}[ll]^{\sigma}.
	}
we find that the projection $\rmpr_{23}$ is dominant and makes $\Gamma$ birationally into a vector bundle over $\PP (V/L)\times \PP (W)$. Hence $\Gamma$
is birational to a succession of vector bundles over $\PP (W)$ or has a \emph{ruled structure} over $\PP (W)$. Since $\bar{G}$ acts generically freely
on $\PP (W)$, the generic fibres of $\varphi$ and $\bar{\varphi}$ can be identified and we can pull back this ruled structure via $\sigma$ (possibly
replacing $\sigma$ by a suitable translate).  Hence  
$\PP (V) /\bar{G}$ is birational to $\PP (W) /\bar{G} \times \PP^N$ with $N = \dim \PP (V) - \dim \PP (W)$. Thus by property (a), $\PP (V)/G$ is
rational.
\end{proof}

In \cite{Shep} essentially this method is used to prove the rationality of the moduli spaces of plane curves of degrees $d\equiv 1$ (mod $9$), $d\ge
19$.\\
It should be noted that covariants are also used in the proof of rationality for $\PP (\mathrm{Sym}^d \CC^2)/\mathrm{SL}_2 (\CC )$, the moduli space
of $d$ points in $\PP^1$, if $d$ is odd and sufficiently large (see \cite{Kat83}).

\subsection[$G$-bundles and configuration spaces of points]{Zero loci of sections in $G$-bundles and
configuration spaces of points}

The technique exposed below was explained to me by P. Katsylo whom I thank for his explanations. The proof is an immediate application of Lemma
\ref{lNoNameLemma}.

\begin{theorem}\xlabel{tZeroLociSections}
Let $G$ be a linear algebraic group and 
let $\mathcal{E}$ be a rank $n$ $G$-vector bundle over a smooth projective $G$-variety $X$ of the same dimension $\dim X= n$. Suppose that $\mathcal{E}$
is spanned by its global sections $V:=H^0(X, \:\mathcal{E})$. Let $N:= c_n (\mathcal{E})$ be the $n$-th Chern class of $\mathcal{E}$. Suppose that the
rational map
\begin{align*}
\alpha\, :\, V &\dasharrow X^{(N)} = \left( \prod_{i=1}^N X \right) /\mathfrak{S}_N\\
s &\mapsto Z(s)
\end{align*}
assigning to a general global section of $\mathcal{E}$ its zeroes, is dominant (thus $\dim H^0 (\mathcal{E}) \ge \dim X \cdot c_n (\mathcal{E})$). If
the action of
$G$ on the symmetric product
$X^{(N)}$ is almost free, then
$V/G$ is birational to
$(X^{(N)}/G)\times
\CC^d$ where
$d =\dim V - N\cdot \dim X$.
\end{theorem}

This result can be applied in two ways: if we know stable rationality of level $\le \dim V -N\cdot \dim X$ of $X^{(N)}/G$, the configuration space of
$N$ unordered points in
$X$, then we can prove rationality of $V/G$. On the other hand, if rationality of $V/G$ is already known, stable rationality of $X^{(N)}/G$ follows.\\
As an example, we consider the space
\[
(\PP^2)^{(7)} / \mathrm{SL}_3 (\CC ) \, ,
\]
the configuration space of $7$ points in $\PP^2$. Rationality of it is proven in the MPI preprint \cite{Kat94}. 

\begin{theorem}\xlabel{tSevenPoints}
The space $(\PP^2)^{(7)} / \mathrm{SL}_3 (\CC )$ is rational.
\end{theorem}

\begin{proof}
If $\mathcal{T}_{\PP^2}$ denotes the tangent bundle of $\PP^2$, then we have $c_2 (\mathcal{T}_{\PP^2}(1))=7$, a general global section of
$\mathcal{T}_{\PP^2}(1)$ has as zero locus seven points in $\PP^2$, and the map
\[
H^0 (\PP^2, \mathcal{T}_{\PP^2}(1)) \dasharrow (\PP^2)^{(7)}
\]
is dominant. Moreover, since $\mathcal{T}_{\PP^2}(1) \simeq \mathcal{R}^{\vee}(2)$, where $\mathcal{R}$ is the tautological subbundle on $\PP^2$ (viewed
as the Grassmannian of $2$-dimensional subspaces in a three-dimensional vector space), we have by the theorem of Borel-Bott-Weil
\[
H^0 (\PP^2, \mathcal{T}(1)) \simeq V(1, \: 2)
\]
as $\mathrm{SL}_3 (\CC )$-representations. Since $\dim V(1, \: 2) = 15$ the map
\[
\PP (H^0 (\PP^2, \mathcal{T}_{\PP^2}(1)))  \dasharrow (\PP^2)^{(7)}
\]
is birational, and
\[
(\PP^2)^{(7)} / \mathrm{SL}_3 (\CC ) \simeq \PP (V (1, \: 2)) /\mathrm{SL}_3 (\CC )\, .
\]
We prove rationality of the latter quotient by a variant of the double bundle method as follows: consider the $\mathrm{SL}_3 (\CC )$-representation 
\[
V = V(1, \: 2) \oplus (V(0,\: 2) \oplus V(1, \: 0)) \oplus V(1, \: 0) \, .
\]
The three-dimensional torus $T= \CC^{\ast} \times \CC^{\ast} \times \CC^{\ast}$ acts in $V$ via
\[
(t_1, \: t_2, \: t_3) \cdot (f, \: (g_1, \: g_2), \: h) = (t_1 f, \: (t_2g_1, \: t_2g_2), \: t_3 h) \, .
\]
We define two $\mathrm{SL}_3$-equivariant maps
\begin{gather*}
\beta \, :\, V(1, \: 2) \times ( V(0, \: 2) \oplus V(1, \: 0)) \to V(1, \: 1), \\
\psi \, : \, V(1, \: 2) \times V(1, \: 0) \to V(1, \: 0) \, .
\end{gather*}

Recall that $V(a, \: b)$ is the kernel of
\[
\Delta = \sum_{i=1}^3 \frac{\partial }{\partial e_i} \otimes \frac{\partial }{\partial x_i}\, :\, \mathrm{Sym}^a (\CC^3)\otimes \mathrm{Sym}^b
(\CC^3)^{\vee }
\to \mathrm{Sym}^{a-1} (\CC^3)\otimes \mathrm{Sym}^{b-1} (\CC^3)^{\vee }
\]
where $e_i$ and $x_j$ are dual coordinates in $\CC^3$ and $(\CC^3)^{\vee }$. In addition there is an $\mathrm{SL}_3 (\CC )$-equivariant map
\begin{gather*}
\omega \, :\, V(a, \: b)\times V(c, \: d) \to \mathrm{Sym}^{a+c+1} (\CC^3) \otimes \mathrm{Sym}^{b+d-2}(\CC^3)^{\vee } \\
\omega (r, \: s) = \sum_{\sigma \in\mathfrak{S}_3} \mathrm{sgn}(\sigma ) e_{\sigma (1)} \frac{\partial r}{\partial x_{\sigma (2)}} \frac{\partial
s}{\partial x_{\sigma (3)} } \, .
\end{gather*}
Then
\begin{gather*}
\beta (f, \: (g_1, \: g_2)) := \Delta (\omega (f,\:  g_1)) + \Delta (fg_2), \: \psi (f, g_2) := \Delta^2 (fg_2^2)\, 
\end{gather*}
(followed by the suitable equivariant projection if necessary).
Thus $\beta$ is bilinear, $\psi$ is linear in the first and quadratic in the second argument. One sets
\begin{gather*}
X:= \{ (f, \: (g_1, g_2), \: h) \, : \, \beta (f, \: (g_1, g_2)) = 0 \; \mathrm{and} \; h \wedge \psi (f, \: g_2)=0 \}\\
\subset V(1, \: 2) \oplus (V(0,\: 2) \oplus V(1, \: 0)) \oplus V(1, \: 0)
\end{gather*}
which is an $\mathrm{SL}_3 (\CC )$ and $T$-invariant subvariety (note that the condition $h \wedge \psi (f, \: g_2)=0$ means that $h$ and $\psi (f, \:
g_2)$ are linearly dependent in $V(1, \: 0) = \CC^3$). For the special points
\begin{gather*}
F = 3e_2x_1x_3 - 2e_1x_1x_2+6e_3x_3x_2-2e_2x_2^2, \: G_1= x_1x_3-x_2^2, \\
G_2 = 2e_2, \: H = \psi (F, G_2) = -32 e_2
\end{gather*}
one checks that $(F, \: (G_1, G_2 ), \: H) \in X$ and that \begin{gather*} \dim \mathrm{ker}(\beta (F, \cdot )) =1,\; \dim \mathrm{ker} (\beta (\cdot ,
(G_1, G_2)))= 7, \\ \dim (\mathrm{ker} (\beta (\cdot, (G_1, G_2))) \cap \mathrm{ker} (\psi (\cdot , G_2))) =4\, .\end{gather*}
So there is a unique irreducible $\mathrm{SL}_3 (\CC )$ and $T$-invariant component $X_0$ of $X$ passing through $(F, \: (G_1, G_2 ), \: H)$; we
consider the two fibration structures on $X_0$ via the projections
\[
\pi_{1}\, :\, X_0 \to V(1, \: 2), \: \pi_2 \, :\, X_0 \to (V(0, \: 2) \oplus V(1, \: 0)) \oplus V(1, \: 0)\, .
\]
The fibres of $\pi_1$ are generically two-dimensional linear spaces in which the subtorus $\{1 \} \times \CC^{\ast }\times \CC^{\ast }$ still acts via
rescaling. Hence
\[
\PP (V(1, \: 2)) /\mathrm{SL}_3 (\CC ) \simeq X_0 / (\mathrm{SL}_3 (\CC )\times T)\, .
\]
On the other hand, via $\pi_2$, $X_0$ is generically a vector bundle over $(V(0, \: 2) \oplus V(1, \: 0)) \oplus V(1, \: 0)$ of rank $5$: in fact,
 $ \dim \mathrm{ker} (\beta (\cdot ,
(G_1, G_2)))= 7$\\ and $\dim (\mathrm{ker} (\beta (\cdot, (G_1, G_2))) \cap \mathrm{ker} (\psi (\cdot , G_2))) =4$, so that the preimage of the line
$\CC H$ in $V(1, \: 0)$ under $\psi (\cdot , G_2)$ restricted to $\mathrm{ker} (\beta (\cdot ,
(G_1, G_2)))$ will be a $5$-dimensional subspace of the $7$-dimensional subspace $\mathrm{ker} (\beta (\cdot ,
(G_1, G_2)))$ of $V(1, \: 2)$. Thus
\begin{gather*}
\PP (V(1, \: 2)) /\mathrm{SL}_3 (\CC ) \simeq \\ \left[ ((V(0, \: 2) \oplus V(1, \: 0)) \oplus V(1, \: 0)) /(\mathrm{SL}_3 (\CC ) \times \CC^{\ast
}\times
\CC^{\ast }) \right] \times \CC^4 \, .
\end{gather*}
But $((V(0, \: 2) \oplus V(1, \: 0)) \oplus V(1, \: 0)) /(\mathrm{SL}_3 (\CC ) \times \CC^{\ast
}\times
\CC^{\ast })$ has dimension $2$. But a unirational surface is rational by Castelnuovo's solution of the L\"uroth problem for surfaces.
\end{proof}

We conclude by remarking that Theorem \ref{tZeroLociSections} makes it obvious, in view of the theorem of Borel-Bott-Weil, that there is an intimate
connection of the rationality problem \ref{problemrat} for a reductive group $G$ and the problem of stable rationality/rationality of configuration spaces of
(unordered) points in generalized flag varieties $G/P$. Since Chern classes of homogeneous bundles $\mathcal{E}$ on $G/P$ arising from a representation
$\varrho\, :\, P
\to \mathrm{Aut} (W)$ can easily be calculated via the splitting principle in terms of the weights of $\varrho$, and $H^0 (G/P, \: \mathcal{E})$ is also
quickly determined by Borel-Bott-Weil, it should not be too difficult to test the range of applicability of Theorem \ref{tZeroLociSections}, but this
remains to be done.

\section[Recent results]{An overview of some recent results on moduli of plane curves}

In \cite{BvB08-1} and \cite{rationalityScript} we prove

\begin{theorem}\xlabel{tBvBmain}
The moduli space $\PP (\mathrm{Sym}^d (\CC^3)^{\vee })/\mathrm{SL}_3 (\CC )$ is rational for all sufficiently large $d$. More precisely it is rational for 
\begin{itemize}
\item[(1)]
$d\equiv 0$ (mod $3$), $d \ge 210$,
\item[(2)]
$d\equiv 1$ (mod $3$), $d \ge 37$,
\item[(3)]
$d\equiv 2$ (mod $3$), $d \ge 65$.
\end{itemize}
\end{theorem}

Part (1) is proven in \cite{Kat89}. For parts (2) and (3) we use the method of covariants as described in Theorem \ref{tCovariants} with the following data: 
$G$ is $\mathrm{SL}_3 (\CC )$ throughout.
\begin{itemize}
\item
For $d= 3n+1$, $n\in\mathbb{N}$, and $V= V(0, \: d)= \mathrm{Sym}^d (\CC^3)^{\vee }$, we take $W = V(0,\: 4)$ and produce covariants
\[
S_d\, :\, V(0,\: d)  \to V(0, \: 4)
\]
of degree $4$. We show that property (b) of Theorem \ref{tCovariants} holds for the space
\[
L_S = x_1^{2n+3} \cdot \CC[x_1,x_2,x_3]_{n-2} \subset V(0, \: d)\, .
\]
Moreover, $\PP (V(0,\: 4))/G$ is stably rational of level $8$. So it suffices to check property (c). 
\item
For $d= 3n+2$, $n\in\mathbb{N}$, and $V= V(0, \: d)= \mathrm{Sym}^d (\CC^3)^{\vee }$, we take $W = V(0,\: 8)$ and produce covariants
\[
T_d\, :\, V(0,\: d)  \to V(0, \: 8)
\]
again of degree $4$. In this case, property (b) of Theorem \ref{tCovariants} can be shown to be true for the subspace
\[
L_T = x_1^{2n+5} \cdot \CC[x_1,x_2,x_3]_{n-3} \subset V(0, \: d)\, .
\]
$\PP (V(0,\: 8))/G$ is stably rational of level $8$, too, hence again everything comes down to checking 
property (c) of Theorem \ref{tCovariants}. 
\end{itemize}

The genericity property (c) of Theorem \ref{tCovariants} is checked in both cases by writing elements in a fibre of $\pi_{L_S}$ (or $\pi_{L_T}$) as
a sum of powers of linear forms, whence $S_d$ resp. $T_d$, being written in Aronhold-Clebsch symbolical notation (see \cite{G-Y}), can be evaluated in a
multilinear fashion. In a final step we reduce modulo a suitable prime $p$ and use upper-semicontinuity over $\mathrm{Spec}(\ZZ )$ to check (c) for all
sufficiently large $d$ not divisible by $3$. For $d=3n+1$ (resp. $d=3n+2$), the calculation to check (c) with the choice of $p=11$ (resp. $p=19$) depends
only on
$n$ modulo
$110$ (resp. $342$), so becomes a finite problem. This is solved computationally in \cite{rationalityScript}.\\
For the details of the argument we refer to \cite{BvB08-1}. The crucial step is to extract data which is sufficient to imply (c) and shows periodic
behaviour over $\mathbb{F}_p$.\\
In \cite{BvB09-1}, also using results from \cite{BvB08-2}, we prove the following result which summarizes our present knowledge on the rationality
properties of the moduli space $C(d)$ of plane curves of degree $d$.

\begin{theorem}\xlabel{tComprehensive}
The moduli space $C(d)$ of plane curves of degree $d$ is rational except possibly for one of the values in the following list: \[ d= 6, \: 7, \: 8, \:
11, \: 12,\:  14,\: 15,\: 16,\: 18,\: 20,\: 23,\: 24,\: 26,\: 32,\: 48 \, . \]
\end{theorem}

The proof is computational and uses the double bundle method for $d$ divisible by $3$, and the method of covariants for $d$ not divisible by $3$. Moreover
for some special values of $d$, e.g. $d=27$, other tricks are employed. To make the problem amenable to the capabilities of present day computers, we had to
use certain new algorithmic approaches, based, roughly speaking, on writing homogeneous polynomials as sums of powers of linear forms and interpolation.
Details appear in \cite{BvB09-1}.

\backmatter

\end{document}